\documentclass[preprint,9pt]{elsarticle}

\usepackage[left=3.5cm,right=3.5cm,top=2.0cm,bottom=2.5cm]{geometry}
\usepackage{graphicx}
\usepackage{amssymb,amsmath,amsfonts}
\usepackage{bm}
\usepackage{titlesec}
\usepackage{changepage} 

\usepackage{float}      
\usepackage{subcaption} 
\usepackage{booktabs}
\usepackage{siunitx} 
\usepackage{float} 
\usepackage{lineno}
\usepackage{xcolor}
\usepackage{lscape}
\usepackage{multirow}
\usepackage{textcomp}
\usepackage{tikz}
\usepackage{tcolorbox}
\usepackage{hyperref}
\hypersetup{
	colorlinks=true,
	linkcolor=blue!70!black,
	urlcolor=blue!70!black,
	citecolor=blue!70!black
}

\usepackage{framed}

\usepackage{titlesec,chapterbib} 
\usepackage{diagbox}

\usepackage{siunitx}

\usepackage{enumitem} 
\definecolor{OliveGreen}{rgb}{0.33,0.71,0.16}

\usepackage{soul}

\usepackage{amssymb}

\journal{}

\begin{document}

\begin{frontmatter}

\title{Scalable multitask Gaussian processes for complex mechanical systems with functional covariates
}

\author[inst1]{Razak C. Sabi Gninkou}
\author[inst2,inst2b]{Andrés~F.~López-Lopera}
\author[inst3]{Franck Massa}
\author[inst4,inst4b]{Rodolphe~Le~Riche}

\affiliation[inst1]{
	organization={Univ. Polytechnique Hauts-de-France, INSA Hauts-de-France, CERAMATHS},
	city={Valenciennes},
	postcode={F-59313},
	country={France}
}

\affiliation[inst2]{
	organization={Univ. de Montpellier, IMAG, CNRS},
	city={Montpellier},
	postcode={F-34090},
	country={France}
}

\affiliation[inst2b]{
	organization={Inria, LEMON},
	city={Montpellier},
	postcode={F-34095},
	country={France}
}

\affiliation[inst3]{
	organization={Univ. Polytechnique Hauts-de-France, CNRS, UMR 8201 - LAMIH, INSA Hauts-de-France},
	city={Valenciennes},
	postcode={F-59313},
	country={France}
}

\affiliation[inst4]{
	organization={CNRS, LIMOS},
	city={Clermont-Ferrand},
	postcode={F-63178},
	country={France}
}

\affiliation[inst4b]{
	organization={École Nationale Supérieure des Mines de Saint-Étienne},
	city={Saint-Étienne},
	postcode={F-42023},
	country={France}
}

\begin{center}
    \begin{abstract}
\begin{adjustwidth}{1.5cm}{1.5cm}
Functional covariates arise in many scientific and engineering applications when
model inputs take the form of time-dependent or spatially distributed profiles,
such as varying boundary conditions or changing material behaviours.
In addition, new practices in digital simulation require predictions accompanied by confidence intervals.
Models based on  Gaussian processes (GPs) provide principled uncertainty quantification.
However, GPs capable of jointly handling functional covariates and
multiple correlated functional  tasks remain largely under-explored. In this work, we extend the framework of GPs with functional covariates to multitask problems by introducing a fully separable kernel structure that captures dependencies across tasks and functional inputs. By taking advantage of the Kronecker structure of the covariance matrix, the model is made scalable.
The proposed model is validated on a synthetic benchmark and 
applied to a realistic structure, a riveted assembly with functional descriptions of the material behaviour and response forces.
The proposed functional multitask GP significantly improves over single task GPs. For the riveted assembly, it requires less than 100 samples to produce an accurate mean and confidence interval prediction. Despite its larger number of parameters, the multitask GP is computationally easier to learn than its single task pendant.

\end{adjustwidth}
\end{abstract}
\end{center}


\end{frontmatter}

\section{Introduction}
\label{sec:intro}

Gaussian processes (GPs) provide a principled Bayesian framework for regression, uncertainty quantification, and surrogate modeling. By defining probability distributions over functions, GPs offer a flexible way to learn complex nonlinear mappings while providing calibrated uncertainty estimates, which is particularly advantageous when data are limited or expensive to acquire.  As a result, GPs have become a well-established approach for statistical modeling across multiple disciplines,
including machine learning~\citep{Rasmussen2006gp}, geostatistics~\citep{Cressie1993statistics},
environmental modeling~\citep{Shen2017gprfunctional}, computational mechanics~\citep{LiDing2026CMAME,BONNEVILLE2021107624}, and robotics~\citep{Deisenroth2015pilco},
where uncertainty-aware surrogate models are essential for decision-making.

In recent years, the role of GP models in computational mechanics has significantly expanded, 
driven by the growing need for surrogate models as an alternative to costly finite element or multi-physics simulations 
in application such as optimization, inverse modeling, reliability analysis, and real-time monitoring \citep{Saida2023transfer, Ogren2024dispersion,LI2020106816}. 
Several studies have demonstrated that GP-based surrogates can accurately reproduce full-field responses and structural quantities 
with only a limited number of high-fidelity simulations, 
enabling efficient uncertainty quantification and reliability assessment in large-scale mechanical systems~\citep{Marrel2024review}. As an example, a fully Bayesian calibration method has been introduced in~\citep{Higdon2008} for computational models with high-dimensional outputs.   
There, both field observations and a limited number of simulator evaluations are jointly modeled 
through a GP emulator defined on a low-dimensional basis. This dimensionality-reduction step makes Bayesian inference and uncertainty quantification feasible for problems involving functional, spatial, or otherwise high-dimensional outputs. The GP paradigm to computer experiments has been extended in~\citep{Qian2008qualitative} to handle both quantitative and qualitative factors, proposing covariance structures that account for similarity across categorical levels while preserving smoothness in continuous inputs. 
This line of works further broadened the applicability of GP models to complex (nonlinear, costly to evaluate) design spaces that mix physical parameters, material types, and operating conditions, making GPs an increasingly relevant setting in computational mechanics. 

Beyond traditional regression tasks, GPs have also been extended to operator learning and inverse problems, enabling the identification of physical input-output mappings governed by high-dimensional or even infinite-dimensional function spaces~\citep{Mora2025OperatorGP,Semler2023adaptive}. 
These latter developments have made GPs not only as black-box surrogates, but also physics-aware probabilistic models capable of capturing structured dependencies in mechanical and multi-physics systems.

Many  applications involve functional inputs rather than scalar. Typical examples in mechanics include force–time loading paths, pressure pulses, temperature-based material properties, or spatially distributed driving fields. In such cases, each experiment or simulation is characterized by an entire function describing the mechanical driver applied to the system \citep{LopezLopera2022cilamce}. The outputs of interest, in contrast, usually correspond to multivariate time series, such as force or displacement responses measured at several structural locations. These responses often exhibit strong temporal correlations and interdependencies across multiple tasks or sensors. Traditional GP models typically assume finite-dimensional  covariates and scalar outputs, which limits their direct applicability in the functional and multitask context. Consequently, surrogate modeling becomes particularly challenging when the covariates must be treated as functions, while the system responses are correlated and coupled across tasks. This setting requires a coherent probabilistic framework capable of jointly modeling functional similarities in the covariates, dependencies in the outputs, and cross-task correlations.

Several approaches have been proposed to handle functional inputs. Classical methods from functional data analysis~\citep{ramsay2005functional} are typically used to reduce infinite-dimensional covariates by projecting them onto a finite set of basis functions such as B-splines, Fourier series, or principal components analysis (PCA),  thereby enabling conventional regression in a finite-dimensional latent space. More recently, GP models in~\citep{Betancourt2020gp,LopezLopera2022,Saunders2021} have incorporated functional similarities directly into tailored kernel definitions, allowing the GPs to operate over spaces of functions rather than finite-dimensional vectors. However, these formulations are  limited to single-output regression.

Conversely, multitask GP (MTGP) models~\citep{Bonilla2008multi}, 
along with their extensions based on process convolutions~\citep{Alvarez2012kernels} 
or graph-based correlations~\citep{Chen2018graph}, are specifically designed to capture 
dependencies across multiple output channels or tasks. To the best of our knowledge, only a limited number of works have attempted to extend 
GP regression to settings involving functional covariates and multivariate outputs. In particular, functional covariates can be incorporated into GP
models through semi-metric constructions defined on infinite-dimensional spaces
\citep{Wang2017functional}.

In many mechanical systems, task-dependent response trajectories are governed by
shared physical mechanisms and common sources of variability, which induce strong
inter-output correlations. Preserving and modelling these correlations is therefore
crucial to ensure coherent uncertainty propagation and physically consistent
predictions for responses indexed by a continuous variable. Ignoring such dependencies
may lead to fragmented predictions and unreliable uncertainty quantification,
especially when responses are sparsely or unevenly observed across tasks.

To address the aforementioned limitations, we introduce an MTGP framework specifically tailored to mechanical systems with functional covariates and multiple correlated tasks. 
The proposed model relies on a fully separable kernel architecture that independently captures similarities across three dimensions:  the task, the functional covariates and the scalar covariate (e.g., the time). 
For tensor-structured data, this construction induces a Kronecker-product decomposition of the covariance matrix, enabling exact yet scalable inference through structured tensor algebra. The resulting implementation, developed in
\texttt{PyTorch}/\texttt{GPyTorch}~\citep{Gardner2018gpytorch}, exploits GPU acceleration for efficient computation
of the marginal likelihood and posterior predictions. 
The model is validated on both synthetic benchmarks and high-fidelity mechanical simulations involving multiple force–displacement response curves driven by nonlinear functional, local, force–displacement profiles. 
To quantify the benefit of our inter-task model, the MTGP is systematically compared against a single-task GP baseline. 

This work is structured as follows. Section~\ref{sec:background} reviews the background of GP modeling with functional covariates and summarizes the challenges associated with functional data. Section~\ref{sec:mtgp} introduces our MTGP formulation for functional inputs and details the construction of the kernel, the parameter estimation, and the prediction. Section~\ref{sec:exp} presents two numerical experiments. We first validate the proposed MTGP with a synthetic example under ideal Gaussian settings, and then test the model on a real-world mechanical application, while assessing its predictive accuracy, its ability to capture cross-task dependencies, and the computational performance. Section~\ref{sec:conclusion} summarizes the main results and gives an overview of relevant future work.

\section{Background}
\label{sec:background}

\subsection{Gaussian processes with functional covariates}

In a functional regression setting, we consider the problem of learning an unknown
deterministic mapping
\begin{equation}
    \begin{array}{rcll}
        y & : & \mathbf{F}(\mathcal{T}, \mathbb{R})^{d_f} &\to \mathbb{R}, \\[0.3em]
          &   & \bm{\mathcal{F}} &\mapsto y(\bm{\mathcal{F}}),
    \end{array}
    \label{eq:background}
\end{equation}
where $\mathbf{F}(\mathcal{T}, \mathbb{R})$ denotes the space of real-valued functions
defined on a compact domain $\mathcal{T} \subset \mathbb{R}$. The functional covariate $\bm{\mathcal{F}} = (f_1, \dots, f_{d_f}) \in \mathbf{F}(\mathcal{T}, \mathbb{R})^{d_f}$
is a $d_f$-dimensional vector of functions, and $y(\bm{\mathcal{F}}) \in \mathbb{R}$
denotes the associated scalar quantity of interest. 

In a mechanical context, the functional covariate $\bm{\mathcal{F}}$ may represent a set
of functional profiles, such as force-displacement responses, contact laws, or
spatially varying fields defined over the domain $\mathcal{T}$. The associated output
$y(\bm{\mathcal{F}})$ may then correspond to a response extracted from the system behavior, for instance a characteristic force level, displacement or velocity fields, or a
failure-related criterion. This functional regression setting serves as a generic starting point. More elaborate
response structures, involving multiple correlated outputs and richer dependency
patterns, will be introduced in Section~\ref{sec:mtgp}.

When modeling computationally expensive computer codes, $y$ is viewed as a black-box function for which only a small number of input evaluations are feasible. To address this challenge, an alternative probabilistic approach is to model 
$y$ as the realization of a stochastic process $\{ Y(\bm{\mathcal{F}}) \}_{\bm{\mathcal{F}} \in \mathbf{F}(\mathcal{T}, \mathbb{R})^{d_f}}$, often assumed to be a GP. If $Y$ is a GP, then, for any finite collection $\bm{\mathcal{F}}_{1}, \ldots, \bm{\mathcal{F}}_{n} \in \mathbf{F}(\mathcal{T}, \mathbb{R})^{d_f}$ with $\bm{\mathcal{F}}_i = (f_{i,1}, \ldots, f_{i,d_f})$, the associated random vector
$\bm{Y}=\big[ Y(\bm{\mathcal{F}}_{1}), \ldots, Y(\bm{\mathcal{F}}_{n}) \big]^{\top} \in \mathbb{R}^{n}$ is multivariate Gaussian-distributed. In particular, under the standard zero-mean prior assumption, we denote by convention $Y \sim \mathcal{GP}(0, k)$, where $k$ is a kernel function defined as  
\begin{equation*}
	\begin{array}{rcl}
		k : 
		& \mathbf{F}(\mathcal{T}, \mathbb{R})^{d_f} \times 
		\mathbf{F}(\mathcal{T}, \mathbb{R})^{d_f}
		& \longrightarrow \mathbb{R}, \\[0.3em]
		& (\bm{\mathcal{F}}, \bm{\mathcal{F}}') 
		& \longmapsto k(\bm{\mathcal{F}}, \bm{\mathcal{F}}') = \operatorname{Cov}\!\big(Y(\bm{\mathcal{F}}), Y(\bm{\mathcal{F}}')\big).
	\end{array}
\end{equation*}

For prediction purposes, we are interested in computing the conditional distribution of \(Y\) given a training set 
$
\mathcal{D} = \{(\bm{\mathcal{F}}_{i}, y_i)\}_{i=1}^{n},
$
where $\bm{\mathcal{F}}_{i}$ denotes the covariate associated to the \(i\)-th code evaluation $y_i \in \mathbb{R}$. Denote $\bm{y} = \big[y_1, \ldots, y_{n} \big]^{\top}$. According to the properties of Gaussian distributions, the predictive distribution of the response at a new input \( \bm{\mathcal{F}}_\star \in \mathbf{F}(\mathcal{T}, \mathbb{R})^{d_f} \) remains Gaussian:
\begin{equation*}
	Y(\bm{\mathcal{F}}_\star) \mid \bm{y} \sim \mathcal{GP}\big(m(\bm{\mathcal{F}}_\star), v(\bm{\mathcal{F}}_\star)\big),
\end{equation*}
where the posterior mean and posterior variance are given by
\begin{align}
    \begin{aligned}
        m(\bm{\mathcal{F}}_\star) &= \bm{k}(\bm{\mathcal{F}}_\star)^\top \mathbf{K}^{-1} \bm{y}, \\
	    v(\bm{\mathcal{F}}_\star) &= k(\bm{\mathcal{F}}_\star, \bm{\mathcal{F}}_\star) - \bm{k}(\bm{\mathcal{F}}_\star)^\top \mathbf{K}^{-1} \bm{k}(\bm{\mathcal{F}}_\star).
    \end{aligned}
    \label{eq:gp_posterior}
\end{align}
with \( \mathbf{K} \in \mathbb{R}^{n \times n} \)  the covariance matrix with entries 
$[\mathbf{K}]_{1 \leq i, j \leq n} = k(\bm{\mathcal{F}}_{i}, \bm{\mathcal{F}}_{j})$, and $\bm{k}(\bm{\mathcal{F}}_\star) \in \mathbb{R}^{n}$ the cross-covariance vector with entries 
$[\bm{k}(\bm{\mathcal{F}}_\star)]_{1 \leq i \leq n} = k(\bm{\mathcal{F}}_\star, \bm{\mathcal{F}}_{i})$.

The closed‑form expressions in Eq.~\eqref{eq:gp_posterior} are a key advantage of GP models, offering both predictive accuracy and uncertainty quantification for complex systems. While the posterior mean  $m(\bm{\mathcal{F}}_\star)$ provides a point prediction of the response at $\bm{\mathcal{F}}_\star$, the posterior variance $v(\bm{\mathcal{F}}_\star)$ quantifies the uncertainty associated to the prediction. 

We should remark that the definition of the posterior covariance 
\begin{equation*}
    \operatorname{cov}(\bm{\mathcal{F}}_\star,\bm{\mathcal{F}}'_\star)
 = k(\bm{\mathcal{F}}_\star,\bm{\mathcal{F}}'_\star)
- \bm{k}(\bm{\mathcal{F}}_\star)^{\top}\mathbf{K}^{-1}\bm{k}(\bm{\mathcal{F}}'_\star)
\end{equation*}
 is often found in the literature. However, in this paper we focus on the posterior variance as this quantity will later be further developed using more efficient Kronecker-based operations in the MTGP framework.

\subsection{Construction of valid kernel}

A central challenge in GP models with functional covariates lies in designing kernels 
that are valid on infinite-dimensional input spaces. A common practice is to consider 
stationary kernel obtained by composing a scalar-valued function 
\(\psi : \mathbb{R}_+ \to \mathbb{R}\) with a dissimilarity measure between functions. 
This leads to kernels of the form
\begin{equation}
    k(\bm{\mathcal{F}}, \bm{\mathcal{F}}') 
    = \psi\!\left( \| \bm{\mathcal{F}} - \bm{\mathcal{F}}' \| \right),
    \label{eqkernel}
\end{equation}
where \( \| \cdot \| \) denotes a dissimilarity measure on 
\(\mathbf{F}(\mathcal{T}, \mathbb{R})^{d_f}\).  
This measure captures both amplitude and shape variations across the domain, and 
the squared structure of \(\| \bm{\mathcal{F}} - \bm{\mathcal{F}}' \|\) makes 
the formulation in Eq.~\eqref{eqkernel} particularly well-suited for radial stationary kernels.

It is important to emphasize that an arbitrary choice of \(\psi\), together with  \(\| \bm{\mathcal{F}} - \bm{\mathcal{F}}' \|\), does not necessarily 
yield a positive semidefinite (psd) kernel. Radial kernels on general metric spaces are guaranteed to be positive definite
provided that \(\psi(\sqrt{\cdot})\) is a completely monotone function
\citep{Schoenberg1938metric}.  Equivalently, $\psi$ admits the scale-mixture representation
\begin{equation}
	\psi(t) = \int_{0}^{\infty} e^{-\omega t^{2}}\,\mu(d\omega),
	\qquad t \ge 0,
\end{equation}
where $\omega>0$ is a scale parameter controlling the smoothness of the elementary  square exponential (SE) kernel $e^{-\omega t^2}$, and $\mu$ is a positive Borel measure on $[0,\infty)$
governing the mixture of scales. This characterization 
provides the theoretical foundation for constructing isotropic kernels in functional settings. Building on these results, the theory of universal kernels shows that such radial
constructions give rise to kernel families that are not only positive semidefinite
but also dense in the space of continuous functions \citep{Micchelli2006universal}.
 We must note that, while alternative dissimilarity measures, such as the 
\(L^\infty\)-norm or Sobolev-type norms, may be considered in theory, they do not guarantee 
positive semidefiniteness when composed with an arbitrary function \(\psi\). Ensuring this 
property typically requires additional compatibility conditions, which are difficult to verify in 
practice~\citep{Sung2022functional}. 

For these reasons, \(L^2\)-based distances remain the most robust and widely used choice in the 
construction of kernels for functional covariates. A particularly well-suited option is the 
weighted \(L^2\)-norm, defined as
\begin{equation}
    \| \bm{\mathcal{F}} - \bm{\mathcal{F}}' \|_{\bm{\ell}}^2 
    = \sum_{d=1}^{d_f} \frac{ \| f_d - f'_d \|_{L^2(\mathcal{T})}^2 }{ \ell_d^2 },
    \label{eq:distance-l2}
\end{equation}
where $\bm{\ell} = (\ell_1, \ldots, \ell_{d_f})$ with \(\ell_d > 0\)  a length-scale parameter controlling the relative sensitivity to variations 
in the \(d\)-th functional component~\citep{Betancourt2020gp,LopezLopera2022}.

Among the kernels that yield valid constructions, a widely used class is the Matérn family which is defined, using the weighted \(L^2\)-norm in~\eqref{eq:distance-l2}, as
\begin{equation}
    \label{eq:MaternKernelFamily}
    k(\bm{\mathcal{F}}, \bm{\mathcal{F}}')
    = \sigma^2 \frac{2^{1-\nu}}{\Gamma(\nu)}
    \left( \sqrt{2\nu} \, \| \bm{\mathcal{F}} - \bm{\mathcal{F}}' \|_{\bm{\ell}} \right)^{\nu}
    K_{\nu}\!\left( \sqrt{2\nu} \, 
    \| \bm{\mathcal{F}} - \bm{\mathcal{F}}' \|_{\bm{\ell}} \right),
\end{equation}
where \(\nu > 0\) and \(\sigma^2 > 0\) denote the smoothness and variance parameters (respectively), 
and where \(\Gamma\) and \(K_\nu\) are the Gamma function and the modified Bessel function 
of the second kind (resp.). 

Specific choices of the smoothness parameter \(\nu\) lead to closed forms. 
For instance, setting \(\nu = \tfrac{5}{2}\) yields the Matérn-\(\tfrac{5}{2}\) kernel,
\begin{equation}
    k(\bm{\mathcal{F}}, \bm{\mathcal{F}}')
    = \sigma^2 \left( 
        1 + \sqrt{5} \, \| \bm{\mathcal{F}} - \bm{\mathcal{F}}' \|_{\bm{\ell}}
        + \frac{5}{3} \| \bm{\mathcal{F}} - \bm{\mathcal{F}}' \|_{\bm{\ell}}^{2}
      \right)
      \exp\!\left( -\sqrt{5} \, 
      \| \bm{\mathcal{F}} - \bm{\mathcal{F}}' \|_{\bm{\ell}} \right),
    \label{marten52}
\end{equation}
and the choices \(\nu = \tfrac{1}{2}\) and \(\nu \to \infty\) recover the exponential and 
the squared-exponential (Gaussian) kernels, respectively~\citep{Genton2001classes}. Hence, \(\nu\) allows to control the regularity of the associated GP: smaller values correspond to rougher sample paths, whereas larger values yield progressively 
smoother trajectories.

In practical settings, functional covariates are only available in discretized or finite-dimensional representations, which calls for suitable approximations of the \( L^2 \)-norm in Eq.~\eqref{eq:distance-l2} to be able to conduct kernel evaluation and model training. Inspired by the techniques on functional data analysis, it is possible to approximate a functional input by a finite-dimensional subspace spanned by an appropriate basis ~\citep{ramsay2005functional}.

\subsection{Dimensionality reduction for functional covariates}
\label{subsec:approximation-functional}
Let \( \{ \Upsilon_{d,r} \}_{r=1}^{p_d} \) be a finite basis of \( L^2(\mathcal{T}) \), where \( \mathcal{T} \subset \mathbb{R} \) is a compact domain. We can then project each functional input $f_d$, for $d \in \{1, \ldots, d_f\}$, onto this basis as
\begin{equation}
	\Pi(f_d)(u) = \sum_{r=1}^{p_d} \beta_{d,r} \Upsilon_{d,r}(u), \label{eq:projection}
\end{equation}
where the projection coefficients \( \bm{\beta}_d = [\beta_{d,1}, \ldots, \beta_{d,p_d}]^\top \) are estimated by the minimization of the \( L^2 \)-norm of the residuals between $f_d$ and $\Pi(f_d)$. This procedure preserves the geometric structure of \( L^2 \), while enabling finite-dimensional computations~\citep{ramsay2005functional,Betancourt2020gp,LopezLopera2022}. 

Using the projection in Eq.~\eqref{eq:projection} for two functions $f_d$ and $f'_d$, then the squared \(L^2\)-distance can be approximated by
\begin{equation}
	\| f_d - f_d' \|_{L^2(\mathcal{T})}^2 \;\approx\; (\bm{\beta}_d - \bm{\beta}'_d)^\top \bm{\Phi}_d (\bm{\beta}_d - \bm{\beta}'_d),
\end{equation}
where $ \bm{\Phi}_d= \int_{\mathcal{T}} \Upsilon_d(u)\,\Upsilon_d^\top(u)\,du$ with $ 
\Upsilon_d(u) = \big[\Upsilon_{d,1}(u),\dots,\Upsilon_{d,p_d}(u)\big]^\top$. Plugging this approximated quantity into in Eq.~\eqref{eq:distance-l2} gives the general approximation
\begin{equation}
	\| \bm{\mathcal{F}} - \bm{\mathcal{F}}' \|_{\bm{\ell}}^2
	\;\approx\;
	\sum_{d=1}^{d_f} \frac{ (\bm{\beta}_d - \bm{\beta}'_d)^\top \bm{\Phi}_d (\bm{\beta}_d - \bm{\beta}'_d) }{ \ell_d^2}.
	\label{eq:distance-l2-approx}
\end{equation}

Two useful simplifications follow from the structure of the Gram matrix $\bm{\Phi}_d$. If the basis $\{\Upsilon_{d,r}\}_{r=1}^{p_d}$ is orthonormal in $L^2(\mathcal{T})$, then $\bm{\Phi}_d=\mathbf{I}$ and
\begin{equation}
	\| \bm{\mathcal{F}} - \bm{\mathcal{F}}' \|_{\bm{\ell}}^2
	\;\approx\;
	\sum_{d=1}^{d_f}  \sum_{r=1}^{p_d}\frac{\big(\beta_{d,r} - \beta'_{d,r}\big)^2}{\ell_d^2}.
	\label{eq:distance-l2-orthonormal}
\end{equation}
If the basis is orthogonal but not normalized, then
$
\bm{\Phi}_d=\operatorname{diag}(\phi_{d,1},\dots,\phi_{d,p_d})$, with 
$\phi_{d,r}=\int_{\mathcal{T}} \Upsilon^2_{d,r}(u) \,du,
$
and thus
\begin{equation}
	\| \bm{\mathcal{F}} - \bm{\mathcal{F}}' \|_{\bm{\ell}}^2
	\;\approx\;
	\sum_{d=1}^{d_f}\sum_{r=1}^{p_d} \frac{\phi_{d,r}\,\big(\beta_{d,r} - \beta'_{d,r}\big)^2}{\ell_d^2}.
	\label{eq:distance-l2-orthogonal}
\end{equation}

A variety of basis function families are commonly used in functional data
analysis, each tailored to specific structural properties of the signals. For
smooth trajectories, B-spline bases provide local support and flexible control of
smoothness \citep{DeBoor2001,Eilers1996pspline}. For periodic or nearly periodic behavior, Fourier bases are a natural choice \citep{ramsay2005functional}. To represent localized and transient features, wavelet bases provide a multiresolution decomposition that captures both time and scale information \citep{Daubechies1992ten,Mallat1999wavelet}. Beyond fixed bases, empirical/data-driven representations such as PCA  reduces dimensionality by learning low-rank structure directly from the data \citep{Jolliffe2016pca}. Variations in the application of PCA have been considered in the literature where functional curves are first represented through a predefined basis expansion after which PCA  is performed  to the resulting coefficient vectors, as illustrated in the functional and wavelet-based analysis in~\citep{Salvatore2016}.

\section{Multitask Gaussian processes with functional covariates}
\label{sec:mtgp}

According to our mechanical application described in Section~\ref{sec:mechanical_application}, and as in many engineering and physical contexts, system responses are influenced not only by global operating conditions, represented here by functional covariates, but also vary with respect to a physical quantity. This additional covariate may correspond to prescribed displacement, time, frequency or a loading-cycle index (e.g., in fatigue analysis). In all such situations, the evolution of the output along this physical quantity cannot be fully captured by the functional covariates alone. For this reason, we consider a physical covariate in addition to the functional ones. Although this covariate could belong to a higher-dimensional space, we will take it as one-dimensional in accordance with our the mechanical application where it is a scalar displacement \(u \in \mathbb{R}\). However, as the notation is kept general, the extension to $u \in \mathbb{R}^{d_u}$ with $d_u > 1$ is straightforward.

Under this setting, we introduce a  multitask framework for modeling multiple correlated tasks $y_s : \mathbf{F}(\mathcal{T}, \mathbb{R})^{d_f} \times \mathbb{R} \to \mathbb{R}$ for $s \in \{1,\ldots,S\}$. 
The core of our formulation lies in a separable kernel structure obtained through tensor products, that decomposes input correlations into components acting on the functional covariates, the scalar domain, and the task index. For tensor-structured datasets, this construction induces a Kronecker product decomposition of the full covariance matrix, which not only enables expressive modeling but also ensures scalability of the inference and prediction procedures.

\subsection{Modeling with a separable kernel structure}
\label{sc:formulation}

We consider a multitask GP (MTGP) regression framework in which each output depends on a set of functional and scalar covariates. Let $(s,\bm{\mathcal{F}}, u) \in \mathcal{S} \times \mathbf{F}(\mathcal{T}, \mathbb{R})^{d_f} \times \mathbb{R}$,
with $\bm{\mathcal{F}} = (f_1, \dots, f_{d_f})$ and $\mathcal{S} = \{1, \dots, S\}$,
denote respectively the index of the corresponding  task, the
$d_f$-dimensional functional covariate, and the scalar covariate, as defined in Section~\ref{sec:background}.

The collection of all task-specific processes \( \{ Y_s \}_{s \in \mathcal{S}} \) is supposed to define $
\bm{Y}_\mathcal{S} = \{\, Y_s(\bm{\mathcal{F}}, u) \,\}_{s \in \mathcal{S}} \sim \mathcal{GP}(0, k)$, where $\bm{Y}_\mathcal{S}$ is GP-distributed vector-valued-function with kernel
\begin{equation}
	k\left((s,\bm{\mathcal{F}}, u),\, (s',\bm{\mathcal{F}}', u')\right)
	= \operatorname{Cov}\left( Y_s(\bm{\mathcal{F}}, u),\, Y_{s'}(\bm{\mathcal{F}}', u') \right),
	\label{eq:mtgp-kernel}
\end{equation}
defined on the extended input space 
\( \mathcal{S} \times \mathbf{F}(\mathcal{T}, \mathbb{R})^{d_f} \times \mathbb{R} \). 
In this formulation, the kernel \( k \) jointly captures dependencies across the task indexes \((s, s')\), the functional covariates  \((\bm{\mathcal{F}}, \bm{\mathcal{F}}')\), and the scalar covariates \((u, u')\).

A common modeling assumption is that correlations induced by the input covariates and by the task structure can be treated independently through a separable kernel construction
 \citep{Bonilla2008multi}. In our case, this assumption leads to
\begin{equation*}
    k\left((s,\bm{\mathcal{F}}, u),\, (s',\bm{\mathcal{F}}', u')\right) = k_{\mathcal{S}}(s,s')\, k_{f, u}\left((\bm{\mathcal{F}}, u),\, (\bm{\mathcal{F}}', u')\right),
\end{equation*}

A similar assumption can also be adopted to further separate the dependence induced by
the functional covariates from that associated with the scalar covariate, as commonly considered in kernel constructions involving mixed variables~\citep{Betancourt2020gp}. Under this setting, the kernel $k$ can be written as

\begin{equation}
    k\!\left((s,\bm{\mathcal{F}}, u),\, (s',\bm{\mathcal{F}}', u')\right)
    = k_{\mathcal{S}}(s,s')\, 
      k_{f}(\bm{\mathcal{F}},\bm{\mathcal{F}}')\,
      k_{u}(u,u'),
    \label{eq:separable-kernel}
\end{equation}
where \(k_{\mathcal{S}}\) is a psd matrix whose entries represent 
additional parameters of the model, \(k_f\) encodes correlations over the functional 
covariates, and \(k_u\) captures correlations along the scalar covariate. 
In what follows, we assume that \(k_f\) and \(k_u\) are valid kernels. Consequently, since the product of psd kernels preserves the psd condition according to the Schur's theorem, the kernel in Eq.~\eqref{eq:separable-kernel} defines a valid covariance function on
\(\mathcal{S} \times \mathbf{F}(\mathcal{T},\mathbb{R})^{d_f} \times \mathbb{R}\) ~\citep{Rasmussen2006gp,Genton2001classes}.

By factorizing the kernel into multiplicative components associated with the task index, the functional and scalar covariates, it captures interaction effects that cannot be represented by purely additive constructions~\citep{Bonilla2008multi,Boashash2003tfd}. The product form further ensures that similarity is preserved simultaneously across all dimensions, which reflects realistic settings where dependencies arise jointly from multiple sources of variation. 
From a computational perspective, the separable structure induces a Kronecker product form in the covariance matrix whenever the sampling design is tensor-structured, as we will discuss in Section \ref{subsec:gp-hyperparam}. This property, as shown in Section \ref{sc:cholesky-likelihood}, will enable efficient algorithms for matrix–vector products, log-determinant evaluations, and linear solves, thereby reducing memory and time complexity compared to dense approaches. Such scalability makes the resulting MTGP framework particularly well suited to large-scale datasets. 

\subsection{Hyperparameter estimation}
\label{subsec:gp-hyperparam}

According to our mechanical context, observations are collected over a fixed discretization grid 
$ \{ u_j \}_{j=1}^{n_{u}} \subset \mathcal{T} $, 
common to all functional covariates and tasks. Each functional replicate is indexed by 
$i \in \{1, \dots, n_f\}$ and each task by 
$s \in \mathcal{S} $, yielding  responses
$Y_{s,i,j} = Y_s(\bm{\mathcal{F}}_{i}, u_j)$, for all triplets \( (s, i, j) \). Stacking all responses gives a single vector 
\(\bm{Y} \in \mathbb{R}^n\) with \(n = S n_f n_{u}\).  
Under the GP assumption, \(\bm{Y}\) follows a joint Gaussian distribution
$
\mathcal{N}(\mathbf{0}, \mathbf{K}_{\boldsymbol{\theta}})$,
where the covariance entries are given by
\begin{equation}
	\left[\mathbf{K}_{\boldsymbol{\theta}}\right]_{(s,i,j),\, (s',i',j')} 
	= k_{\mathcal{S}} \big(s, s') k_{f} \big(\bm{\mathcal{F}}_{i}, \bm{\mathcal{F}}_{i'}) k_{u} \big(u_{j}, u_{j'}),
	\label{eq:covmatrix}
\end{equation}
with $\boldsymbol{\theta}$ denoting the collection of hyperparameters involved in 
$k_{f}$ and $k_{u}$, as well as the parameters defining the positive semidefinite matrix 
$k_{\mathcal{S}}$. For notational simplicity, we omit the explicit dependence of these latter
quantities on $\boldsymbol{\theta}$.

The tensorized indexing associated to  Eq.~\eqref{eq:covmatrix}  implies that $\mathbf{K}_{\boldsymbol{\theta}}$ can be written as a Kronecker product:
\begin{equation}
	\mathbf{K}_{\boldsymbol{\theta}} = \mathbf{K}_{\mathcal{S}} \otimes \mathbf{K}_f \otimes \mathbf{K}_{u},
	\label{eq:kronecker-cov}
\end{equation}
with \([\mathbf{K}_{\mathcal{S}}]_{s,s'} = k_{\mathcal{S}}(s,s')\), 
\([\mathbf{K}_f]_{i,i'} = k_f(\bm{\mathcal{F}}_{i}, \bm{\mathcal{F}}_{i'})\), and  
\([\mathbf{K}_{u}]_{j,j'} = k_{u}(u_j, u_{j'})\).  
Such a decomposition greatly reduces the cost of matrix–vector products and log-determinant computations (see Section~\ref{sc:cholesky-likelihood}), and facilitates modular, interpretable kernel design~\citep{Gilboa2013scaling}. Although Kronecker products are not commutative, the order of the product in Eq.~\eqref{eq:kronecker-cov} is essentially without effect in that it corresponds to a reordering of the data through the ordering of the indices triplets ($s$, $\bm{\mathcal{F}}$, $u$).

By considering the Mat\'ern kernel structure in Eq.~\eqref{eq:MaternKernelFamily} for $k_f$ and adopting an analogous kernel construction for $k_u$ based on the Euclidean norm $\|u - u'\|_{\ell_u}$, then we have the set of hyperparameters
\begin{equation}
	\boldsymbol{\theta} = 
	\left( \mathbf{K}_{\mathcal{S}}, \sigma^2, \boldsymbol{\ell}_f, \ell_{u} \right), \label{parameter}
\end{equation}
where \(\mathbf{K}_{\mathcal{S}} \in \mathbb{R}^{S \times S}\) encodes inter-task correlations,  
\(\boldsymbol{\ell}_f = (\ell_1, \dots, \ell_{d_f}) \in \mathbb{R}^{d_f}_+\) and \( \ell_{u} > 0\) are the length-scale parameters of the functional and scalar kernels (resp.), and \(\sigma^2 > 0\) is a global variance parameter. We note that separate variance parameters for the kernels $k_f$ and $k_u$ are not considered here, as this would lead to an identifiability issue: the overall variance of $k$ would then be given by the product $\sigma^{2} = \sigma_f^{2}\sigma_u^{2}$.

Although the Eq.~\eqref{parameter} suggests that \(\mathbf{K}_{\mathcal{S}}\) could be estimated as a 
full symmetric psd matrix, this is rarely done in practice. A full-rank parametrization 
involves \(S(S+1)/2\) free parameters and often leads to poor identifiability and numerical 
instabilities, especially as \(S\) increases. An alternative in MTGP modeling is to factorize \(\mathbf{K}_{\mathcal{S}}\) via Cholesky decomposition  
\begin{equation}
\mathbf{K}_{\mathcal{S}} = \mathbf{L}_{\mathcal{S}}\mathbf{L}_{\mathcal{S}}^\top,
\qquad \mathbf{L}_{\mathcal{S}}\in\mathbb{R}^{S\times S},
\label{eq:chol_tasks}
\end{equation}
which guarantees positive definiteness (up to numerical tolerances), and can improve
optimization stability. 
We note that it is also possible to consider a low-rank
factorization $\mathbf{K}_{\mathcal{S}}=\mathbf{C}_{\mathcal{S}}\mathbf{C}_{\mathcal{S}}^\top$ with $\mathbf{C}_{\mathcal{S}}\in\mathbb{R}^{R\times S}$ and $R\ll S$,
to reduce the amount of parameters that need to be estimated. However this later option is not further considered here as it does not allow further simplifications in the computation of the determinant and inversion of $\mathbf{K}_{\mathcal{S}}$, both quantities required in the likelihood evaluation (See Section~\ref{sc:cholesky-likelihood}).

We now turn to the estimation of 
\(\boldsymbol{\theta}\), which is performed by minimizing the negative log-marginal likelihood associated with the MTGP model given by 
\begin{equation}
    \mathcal{L}(\boldsymbol{\theta})
    =  \frac{1}{2} \log |\mathbf{K}_{\boldsymbol{\theta}}| +\frac{1}{2} \bm{y}^\top \mathbf{K}_{\boldsymbol{\theta}}^{-1} \bm{y}+ \frac{n}{2} \log(2\pi),
    \label{likelihood}
\end{equation}
with the observation vector \(\bm{y} = [y_1,\ldots,y_{n}]^{\top}\). This quantity is typically minimized using gradient-based optimizers such as L-BFGS~ or Adam~\citep{Liu1989limited,Kingma2015adam}.  From~\eqref{likelihood}, we note that evaluating the $\mathcal{L}$ requires repeated computations of 
\(\mathbf{K}_{\boldsymbol{\theta}}^{-1}\) and 
\(\log |\mathbf{K}_{\boldsymbol{\theta}}|\).  
For large multitask problems, these operations become prohibitively expensive unless the Kronecker structure of 
\(\mathbf{K}_{\boldsymbol{\theta}}\) is fully exploited.

\subsection{Efficient Kronecker-based inference}
\label{sc:cholesky-likelihood}
A key property of Kronecker products is that the Cholesky decomposition factorizes across 
tensor components.  
Specifically, if $\mathbf K_{d}=\mathbf L_{d}\mathbf L_{d}^\top$, $d\in\{\mathcal S,f,u\}$, then the global covariance matrix admits the Cholesky factorization $ \mathbf{K}_{\boldsymbol{\theta}} 
    = \mathbf{L}\,\mathbf{L}^\top$,
   $ \mathbf{L} = 
    \mathbf{L}_{\mathcal{S}} \otimes \mathbf{L}_{f} \otimes \mathbf{L}_{u}
    \in \mathbb{R}^{n \times n}$,
with \(n = S n_f n_{u}\).  
This avoids computing the Cholesky decomposition of the full matrix 
$\mathbf{K}_{\boldsymbol{\theta}}$, replacing it with the much cheaper decompositions of 
$\mathbf{K}_{\mathcal{S}}$, $\mathbf{K}_{f}$, and $\mathbf{K}_{u}$.

Using the Cholesky factorization, then Eq.~\eqref{likelihood} can be rewritten as
\begin{equation}
    \mathcal{L}(\boldsymbol{\theta})
    = \log |\mathbf{L}|
    + \frac{1}{2} \|\boldsymbol{\alpha}\|^{2}
    + \frac{n}{2} \log(2\pi),
    \label{eq:nlml_cholesky}
\end{equation}
where \(\mathbf{L} \boldsymbol{\alpha} = \bm{y}\) is a linear system that can be solved by forward substitution.  The key computational challenge therefore lies in evaluating both 
\(\log |\mathbf{L}|\) and \(\boldsymbol{\alpha}\) efficiently, with the latter being the 
most demanding operation.

For the computation of \(\log |\mathbf{L}|\) (see~\ref{appendix:Kronecker:logL}), 
the separable structure yields the closed-form expression
\begin{equation}
    \log |\mathbf{L}|
    = (n_f n_u)\sum_{i=1}^{S} \log (\mathbf{L}_{\mathcal{S}})_{ii}
    + (S n_u)\sum_{i=1}^{n_f} \log (\mathbf{L}_{f})_{ii}
    + (S n_f)\sum_{i=1}^{n_u} \log (\mathbf{L}_{u})_{ii}.
\end{equation}
This formulation avoids computing \(|\mathbf{L}|\) by first forming the full Cholesky factor
\(\mathbf{L} \in \mathbb{R}^{n \times n}\), a step that can be time-consuming. Instead, the computation reduces to the sum of the logarithms of the diagonal entries of the Cholesky factors \(\mathbf{L}_{\mathcal{S}}\), \(\mathbf{L}_{f}\), and \(\mathbf{L}_{u}\), calculations that can be achieved more efficiently.
 
We now turn to the computation of $\boldsymbol{\alpha}$ which is given by solving the linear system \(\mathbf{L} \boldsymbol{\alpha} = \bm{y}\)
with $\mathbf{L} = \mathbf{L}_{\mathcal{S}} \otimes \mathbf{L}_{f} \otimes \mathbf{L}_{u}$. Using the inversion rule for Kronecker products, we obtain
\begin{equation*}
    \boldsymbol{\alpha}
    = (\mathbf{L}_{\mathcal{S}}^{-1} \otimes \mathbf{L}_{f}^{-1} \otimes \mathbf{L}_{u}^{-1})\,\bm{y}.
\end{equation*}

To this end, the vector \(\bm{y} \in \mathbb{R}^n\) is reshaped into the tensor
\(\mathcal{Y} \in \mathbb{R}^{S \times n_f \times n_u}\), where the modes correspond to 
tasks, functional replicates, and temporal grid points.  
Using the identity relating Kronecker products and vectorization 
(see~\ref{Kronecker--vec}), we obtain
\begin{equation}
    \boldsymbol{\alpha}
    = \operatorname{vec}\!\left(
        \mathcal{Y}
        \times_{u} \mathbf{L}_{u}^{-1}
        \times_{f} \mathbf{L}_{f}^{-1}
        \times_{\mathcal{S}} \mathbf{L}_{\mathcal{S}}^{-1}
    \right),
    \label{eq:alpha}
\end{equation}
where \(\times_d\) denotes the mode-wise multiplication along dimension \(d\). These mode-wise multiplications can be achieved by solving the following linear systems:
\begin{equation}
    \label{eq:modewise-prod}
    \begin{aligned}
        \mathbf{L}_{u} (\mathcal{A}^{(1)})_{s,i,:} &= (\mathcal{Y})_{s,i,:}, 
        && \forall (s,i),\\[0.3em]
        \mathbf{L}_{f} (\mathcal{A}^{(2)})_{s,:,j} &= (\mathcal{A}^{(1)})_{s,:,j}, 
        && \forall (s,j),\\[0.3em]
        \mathbf{L}_{\mathcal S} (\mathcal{A}^{(3)})_{:,i,j} &= (\mathcal{A}^{(2)})_{:,i,j}, 
        && \forall (i,j),\\[0.6em]
        \boldsymbol{\alpha} &= \operatorname{vec}\!\big(\mathcal{A}^{(3)}\big), &&
    \end{aligned}
\end{equation}
where $\mathcal{A}^{(1)}, \mathcal{A}^{(2)}, \mathcal{A}^{(3)}$ are intermediate tensors produced after each mode-wise triangular solve. More precisely, the computation proceeds through a sequence of mode-wise triangular solves.
First, the lower triangular system
$\mathbf{L}_{u} (\mathcal{A}^{(1)})_{s,i,:} = (\mathcal{Y})_{s,i,:}$
is solved along the scalar mode for all $(s,i)$, yielding the tensor
$\mathcal{A}^{(1)}$.
Similarly, the systems
$\mathbf{L}_{f} (\mathcal{A}^{(2)})_{s,:,j} = (\mathcal{A}^{(1)})_{s,:,j}$ and $\mathbf{L}_{\mathcal{S}} (\mathcal{A}^{(3)})_{:,i,j} = (\mathcal{A}^{(2)})_{:,i,j}$ are solved along the functional mode for all $(s,j)$ to have $\mathcal{A}^{(2)}$ and along the task mode for all $(i,j)$ to obtain $\mathcal{A}^{(3)}$.

Note from Eq.~\eqref{eq:modewise-prod} that the explicit computation of $\mathbf{L}^{-1}$ is avoided, and that the overall complexity is reduced to solving linear systems involving Cholesky factors of smaller sizes. Consequently, this tensor-aware computation reduces the naive cubic complexity $\mathcal{O} \big((S\,n_f\,n_{u})^3\big)$ to $\mathcal{O} \big(
        S\,n_f\,n_{u}^2
        + S\,n_{u}\,n_f^2
        + n_f\,n_{u}\,S^2
    \big) + 
    \mathcal{O} \big(S^3 + n_f^3 + n_{u}^3\big)$,  where the first term corresponds to the three mode-wise triangular solves,  
and the second to the small Cholesky factorizations of 
\(\mathbf{K}_{\mathcal S}\), \(\mathbf{K}_f\), and \(\mathbf{K}_u\). Regarding memory usage, the requirement decreases from  
$\mathcal{O} \big((S\,n_f\,n_u)^2\big)$ to $\mathcal{O} \big(S^2 + n_f^2 + n_u^2\big)$,
since only the kernel blocks and their Cholesky factors need to be stored.

\subsection{Posterior prediction}
\label{app:posterior-derivation}

Let $\bm{y} \in \mathbb{R}^n$ be the observation vector, with $n = S n_f n_{u}$ as defined in Section \ref{subsec:gp-hyperparam}. 
For a new test input $(\bm{\mathcal{F}}_\star,u_\star)\in \mathbf{F}(\mathcal{T}, \mathbb{R})^{d_f} \times \mathbb{R}$ and a task index $s \in \{1,\dots,S\}$, 
the joint prior distribution of the training and test outputs is Gaussian:
\begin{equation*}
	\begin{bmatrix}
		\bm{y} \\[0.2em] Y_s(\bm{\mathcal{F}}_\star,u_\star)
	\end{bmatrix}
	\sim \mathcal{N}\!\left(
	\mathbf{0},\;
	\begin{pmatrix}
		\mathbf{K} & \bm{k}_\star^{(s)} \\
		\bm{k}_\star^{(s)\top} & k_s(\bm{\mathcal{F}}_\star,u_\star)
	\end{pmatrix}
	\right),
\end{equation*}
where $\bm{k}_\star^{(s)}=\left[k\left((s,\bm{\mathcal{F}}_\star,u_\star),(1,\bm{\mathcal{F}}_1, u_1)\right), \ldots ,k\left((s,\bm{\mathcal{F}}_\star, u_\star),(s,\bm{\mathcal{F}}_{n_f}, u_{n_u})\right) \right]^\top \in \mathbb{R}^{n}$ is the cross-covariance vector between the training set and the test point, $\mathbf{K} \in \mathbb{R}^{n\times n}$ is the training covariance matrix,
and $k_s(\bm{\mathcal{F}}_\star,u_\star) = k\left((s,\bm{\mathcal{F}}_\star,u_\star),(s,\bm{\mathcal{F}}_\star,u_\star)\right) = \sigma^2$ is the prior variance at this location. For notation simplicity,  we intentionally dropped the dependency on the hyperparameters $\boldsymbol{\theta}$ from the indices of the quantities implying the evaluation of the kernel $k$.

By the Gaussian properties, 
the posterior distribution of $Y_s(\bm{\mathcal{F}}_\star,u_\star)$ 
is also Gaussian:
\begin{equation*}
	Y_s(\bm{\mathcal{F}}_\star,u_\star)\mid \bm{y} 
	\;\sim\; 
	\mathcal{GP}\!\big(
	m_s(\bm{\mathcal{F}}_\star,u_\star),\;
	v_s(\bm{\mathcal{F}}_\star,u_\star)
	\big),
\end{equation*}
with posterior mean and variance are given by
\begin{align}
	\begin{aligned}
	    m_s(\bm{\mathcal{F}}_\star,u_\star) 
		&= \bm{k}_\star^{(s)\top} \mathbf{K} ^{-1} \bm{y},\\
		v_s(\bm{\mathcal{F}}_\star,u_\star) 
		&= \sigma^2 -\bm{k}_\star^{(s)\top} \mathbf{K}^{-1} \bm{k}_\star^{(s)}.
	\end{aligned}
\end{align}

The direct computation of $m_s(\bm{\mathcal{F}}_\star,u_\star)$ and 
$v_s(\bm{\mathcal{F}}_\star,u_\star)$ would normally require the explicit inversion of the 
covariance matrix $\mathbf{K}$, which is computationally 
prohibitive for large-scale multitask settings. 
Instead, we can again exploit the Cholesky factorization
$\mathbf{K} = \mathbf{L}\mathbf{L}^\top$ and reuse the vector 
$\boldsymbol{\alpha}$ already computed during likelihood 
optimization. The key efficiency gain arises from the Kronecker-separable structure of the kernel, 
which allows both $\bm{k}_\star^{(s)}$ and $\mathbf{L}$ to be decomposed into mode-wise components. For the former term corresponding to the cross-covariance vector, we have
$$\bm{k}_\star^{(s)} 
= [\bm{k}_{\mathcal{S}}]_s \otimes \bm{k}_f(\bm{\mathcal{F}}_\star) \otimes \bm{k}_{u}(u_\star),
$$
with 
$[\bm{k}_{\mathcal{S}}]_s= \big[k_{\mathcal{S}}(s,1),\,\dots,\,k_{\mathcal{S}}(s,S)\big]^\top$, 
	$\bm{k}_f(\bm{\mathcal{F}}_\star)= \big[k_f(\bm{\mathcal{F}}_\star,\bm{\mathcal{F}}_1),\,\dots,\,k_f(\bm{\mathcal{F}}_\star,\bm{\mathcal{F}}_{n_f})\big]^\top$, and
	$\bm{k}_{u}(u_\star) =  \big[k_{u}(u_\star, u_1),\,\dots, \, k_{u}(u_\star, u_{n_{u}})\big]^\top$.   As $\mathbf{L}^{-1} = \mathbf{L}_{\mathcal{S}}^{-1} \otimes \mathbf{L}_f^{-1} \otimes \mathbf{L}_{u}^{-1}$, combining both Kronecker decompositions yields
$$
\mathbf{L}^{-1}\bm{k}_\star^{(s)} 
= \big(\mathbf{L}_{\mathcal{S}}^{-1}[\bm{k}_{\mathcal{S}}]_s\big)\;\otimes\;
\big(\mathbf{L}_f^{-1}\bm{k}_f(\bm{\mathcal{F}}_\star)\big)\;\otimes\;
\big(\mathbf{L}_{u}^{-1}\bm{k}_{u}(u_\star)\big)
= \bm{\zeta}_{\mathcal{S}}\otimes  \bm{\zeta}_f \otimes  \bm{\zeta}_{u}, 
$$
where the linear systems $\mathbf{L}_{\mathcal{S}}  \bm{\zeta}_{\mathcal{S}} = [\bm{k}_{\mathcal{S}}]_s$, $\mathbf{L}_f  \bm{\zeta}_f = \bm{k}_f(\bm{\mathcal{F}}_\star)$ and $\mathbf{L}_{u}  \bm{\zeta}_{u} = \bm{k}_{u}(u_\star)$ can be solved by forward substitution. Using the above result and Eq.~\eqref{eq:alpha}, the posterior mean is then given by
\begin{equation}
	m_s(\bm{\mathcal{F}}_\star,u_\star) 
	= (\mathbf{L}^{-1}\bm{k}_\star^{(s)})^\top \boldsymbol{\alpha}
	= \big(\bm{\zeta}_{\mathcal{S}}^\top \otimes \bm{\zeta}_f^\top \otimes \bm{\zeta}_{u}^\top\big)\,\boldsymbol{\alpha},
	\label{eq:post-mean-final}
\end{equation}
and the posterior variance follows
\begin{align}
    v_s(\bm{\mathcal{F}}_\star,u_\star) 
    = \sigma^2 - (L^{-1} \bm{k}_\star^{(s)} )^\top L^{-1}\bm{k}_\star^{(s)}
    =  \sigma^2 - (\bm{\zeta}_{\mathcal{S}}^{\top}\bm{\zeta}_{\mathcal{S}})
    (\bm{\zeta}_f^{\top}\bm{\zeta}_f)
    (\bm{\zeta}_u^{\top}\bm{\zeta}_u).
    \label{eq:post-var-final}
\end{align}
where we used the Kronecker property $(\bm{\zeta}_{\mathcal{S}} \otimes \bm{\zeta}_f \otimes \bm{\zeta}_u)^{\!\top}
    (\bm{\zeta}_{\mathcal{S}} \otimes \bm{\zeta}_f \otimes \bm{\zeta}_u) = (\bm{\zeta}_{\mathcal{S}}^{\top}\bm{\zeta}_{\mathcal{S}})
    (\bm{\zeta}_f^{\top}\bm{\zeta}_f)
    (\bm{\zeta}_u^{\top}\bm{\zeta}_u)$. 

Note from Eq.~\eqref{eq:post-mean-final} and Eq.~\eqref{eq:post-var-final} that both  quantities thus rely on mode-wise triangular solves and simple matrix–vector contractions, and therefore, ensuring scalability while preserving numerical stability.

\section{Numerical experiments}
\label{sec:exp}

\subsection{Experimental setup}
\label{subsec:experim_setup}
Implementations are based on the \texttt{GPyTorch} library
\citep{Gardner2018gpytorch}, and all experiments were conducted on a
workstation equipped with an Intel(R) Core(TM) Ultra 7 155H processor
(16 physical cores, 22 logical processors, up to 4.8\,GHz) and
30\,GB of RAM, running Ubuntu 24.04 LTS.
The models were implemented in Python~3.12.2 using
PyTorch~2.5.1 and GPyTorch~1.14. The GPyTorch toolbox has been adapted to compute the separable kernel
structure introduced in Eq.~\eqref{eq:separable-kernel}, relying on the
Kronecker-based implementation described in
Section~\ref{sc:cholesky-likelihood}. All the source codes, and the notebook required to partially reproduce the toy example in Section~\ref{subsec:toy-dataset}, are
publicly available at
\url{https://github.com/SABI-GNINKOU/F-MTGP}.

The predictive performance of the model is evaluated on two datasets:  
(i) a synthetic example designed to provide interpretable behavior and controlled variability 
(Section~\ref{subsec:toy-dataset}), and  
(ii) a real dataset originating from a complex mechanical component 
(Section~\ref{sec:mechanical_application}). For the synthetic  example, a single functional encoding strategy is adopted.
Specifically, the functional inputs are projected onto a low-dimensional latent space
using  PCA, which provides a simple and interpretable
baseline in a controlled setting. In contrast, the mechanical application involves functional inputs with more diverse
profiles in terms of amplitude, smoothness, and frequency content. To account for this variability, several functional dimensionality reduction strategies
are considered and compared. These include  PCA, functional basis expansions based on a Haar wavelet basis
(level~4),  B-splines, as well as hybrid approaches in which PCA is applied to the
corresponding basis coefficients, named here as ``Wavelet + PCA'' and ``B-spline + PCA''. A brief description of the previous encoding strategies is given in~\ref{app:encodings}. In all cases, the functional inputs are mapped to a latent space of fixed dimension
$d_{\text{proj}} = 6$, ensuring a fair comparison between representations of comparable
complexity.
For direct PCA, the six dimensions capture at least $99.9\%$ of the total variance
across the three functional components.
For wavelet- and B-spline-based encodings, the six most energetic coefficients are retained. In the hybrid variants a whitened PCA with six latent directions also explains
approximately $99.9\%$ of the total energy.

Regardless of the chosen functional encoding, all models are trained under the same optimization protocol. Regarding the MTGP, kernel operations rely on Kronecker-structured linear algebra, enabling efficient evaluation of the marginal likelihood and its gradients. In all experiments, the hyperparameters are estimated by maximizing the marginal
log-likelihood using the Adam optimizer \citep{Kingma2015adam}, with a learning rate
$\eta = 2 \times 10^{-2}$, momentum parameters $(\beta_1,\beta_2) = (0.98, 0.999)$,
and a weight decay of $10^{-5}$.
Gradient norms are clipped to a maximum value of $1.0$ to improve numerical stability
during training \citep{pascanu2013difficulty}, and an early-stopping criterion is applied
when the log-likelihood improvement remains below $10^{-3}$ for $20$ consecutive iterations. The maximum number of optimization iterations is set to $n_{\max} = 5 \times 10^{2}$
for the synthetic Rayleigh-based example and $n_{\max} = 2 \times 10^{4}$ for the
mechanical application.
The larger iteration budget allocated to the mechanical case reflects its higher
dimensionality and increased complexity, which require additional iterations to reach stable likelihood values.

Training the MTGP model for the mechanical application is expected to be more challenging
than for the synthetic data, as the data--model consistency is not guaranteed in the
presence of experimental noise, modeling errors, and complex physical interactions.
To mitigate the risk of convergence to sub-optimal local maxima of the marginal
log-likelihood, a multi-start optimization strategy is therefore employed for the
mechanical case, using $n_{\text{restart}} = 10$ random initializations. Before each restart, kernel hyperparameters are reinitialized within ranges adapted to the
 observed variability of the mechanical responses shown in Figure~\ref{fig:functional_inputs_outputs_meca}:
\[
\ell_f \sim \mathcal{U}([0.5, 20]), 
\qquad 
\ell_u \sim \mathcal{U}([0.005, 0.1]), 
\qquad
\sigma^2 \sim \mathcal{U}([0.5, 2]), 
\qquad
\sigma^2_{\mathrm{noise}} \sim \mathcal{U}([10^{-3},10^{-1}]).
\]
Among all restarts, the model achieving the highest log-likelihood value is retained
for subsequent predictions.

\subsection{Performance criteria}
\label{subsec:metrics}

The predictive performance of the 
MTGP model is evaluated on a set of test evaluation points
$\{(\bm{\mathcal{F}}_{\star,i}, u_{\star})\}_{i=1}^{n_{\mathrm{test}}}$
using two complementary criteria:
the coefficient of determination ($Q^2$)
and the coverage accuracy (CA). For each task \(s \in \{1, \ldots, S \}\), the MTGP predictive distribution at a test pair
\((\bm{\mathcal{F}}_{\star,i}, u_{\star})\) is characterized by its conditional
posterior mean
\(m_s(\bm{\mathcal{F}}_{\star,i},u_{\star})\) and variance
\(v_s(\bm{\mathcal{F}}_{\star,i},u_{\star})\), as given in
Eq.~\eqref{eq:post-mean-final} and Eq.~\eqref{eq:post-var-final}.
These predictions are compared to the corresponding true simulator outputs
\(y_s^{\mathrm{true}}(\bm{\mathcal{F}}_{\star,i}, u_{\star})\).

The determination coefficient is defined as
\begin{equation}
Q^2_s
=
1 -
\frac{
\sum_{i=1}^{n_{\mathrm{test}}}
\left(
y_s^{\mathrm{true}}(\bm{\mathcal{F}}_{\star,i},u_{\star})
-
m_s(\bm{\mathcal{F}}_{\star,i},u_{\star})
\right)^2
}{
\sum_{i=1}^{n_{\mathrm{test}}}
\left(
y_s^{\mathrm{true}}(\bm{\mathcal{F}}_{\star,i}, u_{\star})
-
\overline{y}_s
\right)^2
},
\end{equation}
where \(\overline{y}_s\) denotes the empirical mean of all test  for task~\(s\). This quantity measures the accuracy of the point predictions relative to the intrinsic variability of the data. Values of \(Q^2_s\) close to~1 indicate accurate predictions.

The calibration of the predictive uncertainty is quantified through the coverage accuracy of the 
predictive interval
\[
I_{s}(\bm{\mathcal{F}},u, \delta)
=
\left[
m_s(\bm{\mathcal{F}},u)
-
\delta \sqrt{v_s(\bm{\mathcal{F}},u)};
\;
m_s(\bm{\mathcal{F}}, u)
+
\delta \sqrt{v_s(\bm{\mathcal{F}}, u)}
\right],
\]
defined as
\begin{equation}
\mathrm{CA}_s(\delta)
=
\frac{1}{n_{\mathrm{test}}}
\sum_{i=1}^{n_{\mathrm{test}}}
\mathbf{1}_{\left\{
y_s^{\mathrm{true}}(\bm{\mathcal{F}}_{\star,i}, u_{\star,i})
\in 
I_s(\bm{\mathcal{F}}_{\star,i}, u_{\star}, \delta)
\right\}}.
\end{equation}
This criterion estimates the proportion of test points lying inside the predictive
intervals. In our experiments, we set $\delta=1.96$, which corresponds to a nominal $95\%$
predictive interval under Gaussian assumptions, so that a well-calibrated model should
satisfy $\mathrm{CA}_s(1.96)\approx 0.95$.

\subsection{The Rayleigh-based synthetic benchmark}
\label{subsec:toy-dataset}

\subsubsection{Dataset generation}
We construct a synthetic dataset to evaluate both the predictive accuracy and
the scalability of the MTGP model. To guarantee that the problem is well-posed and avoid complexities related to data-model mismatch \citep{Karvonen2023},
the dataset is constructed from samples of a baseline MTGP with known parameters.  It is then certain that the data can be learned by the model. 

Inspired by the upcoming mechanical application in Section~\ref{sec:mechanical_application}, we consider inputs $\bm{\mathcal{F}} = (f_{1}, f_{2}, f_{3})$ where each $f_i$ follows a Rayleigh-shaped
curve \begin{equation*}
    h_\rho(u) = \frac{u}{\rho^2} \exp\left(-\frac{u^2}{2\rho^2}\right),
\end{equation*}
with 
$u \sim \mathcal U([0,1.5])$ and $\rho \sim \mathcal U([0.05,1])$.
To control the amplitude of the function, we use the rescaled form
\begin{equation*}
    f_i(u) \sim \frac{\alpha \, h_\rho(u)}{\max h_\rho(u)},
\end{equation*}
 with $\alpha \sim \mathcal U([2,4])$. These functions are discretized over 150 points
that is, each functional input
$
\bm{\mathcal{F}}_i = (f_{i,1}, f_{i,2}, f_{i,3})
\in \mathbf{F}(\mathcal{T},\mathbb{R})^{3}
$
consists of three full curves, represented after discretization as
$
f_{i,d} \;\longmapsto\;
\big(f_{i,d}(u_1), \dots, f_{i,d}(u_{150})\big) \in \mathbb{R}^{150},
\ d \in \{1,2,3\}
$. Thus, there are \(d_f = 3\) functional channels and each input
\(\bm{\mathcal{F}}_i\) is a high-dimensional functional object. 
For the subsequent steps the generated inputs are projected onto a 6-dimensional PCA basis.

The baseline MTGP is equipped with the separable kernel structure of Eq.~\eqref{eq:separable-kernel} where $k_f$ is a Matérn-$5/2$ kernel as defined in Eq.~\eqref{marten52} with fixed length-scales $\boldsymbol{\ell}_f=(80,80,80)$. For the scalar kernel, we use the additive structure 
$
    k_u(u,u') = k_{\text{Mat}}(u,u') \;+\; k_{\text{Per}}(u,u')
$ 
where 
\begin{equation*}
    k_{\text{Mat}}(u,u') = 
\left(
1 
+ \frac{\sqrt{5}\,|u-u'|}{\ell_{\text{Mat}}}
+ \frac{5\,|u-u'|^{2}}{3\,\ell_{\text{Mat}}^{2}}\right)
\exp\left(-\frac{\sqrt{5}\,|u-u'|}{\ell_{\text{Mat}}}\right),
\end{equation*}
 and \begin{equation*}
     k_{\text{Per}}(u,u') = 
\exp\left(
 -\,\frac{2}{\ell_{\text{Per}}^{2}}\,
 \sin^{2}\left(\frac{\pi |u-u'|}{p}\right)\right),
 \end{equation*} 
 with parameters $ 
\ell_{\text{Mat}}=1.5$, $
\ell_{\text{Per}}=0.5$,
$p=1$. This choice combines the smooth long-range behavior of a Matérn-\(5/2\) kernel with the oscillatory structure induced by a periodic component.  
Such an additive kernel reproduces both transient and quasi-periodic temporal patterns, partially matching the dynamics observed in the real structural component studied in Section~\ref{sec:mechanical_application}. The baseline MTGP is evaluated at \(n_f = 500\) distinct functional inputs.

\begin{figure}[t!]
	\centering
	\begin{subfigure}[t]{0.325\linewidth}
		\centering
		\includegraphics[width=\linewidth]{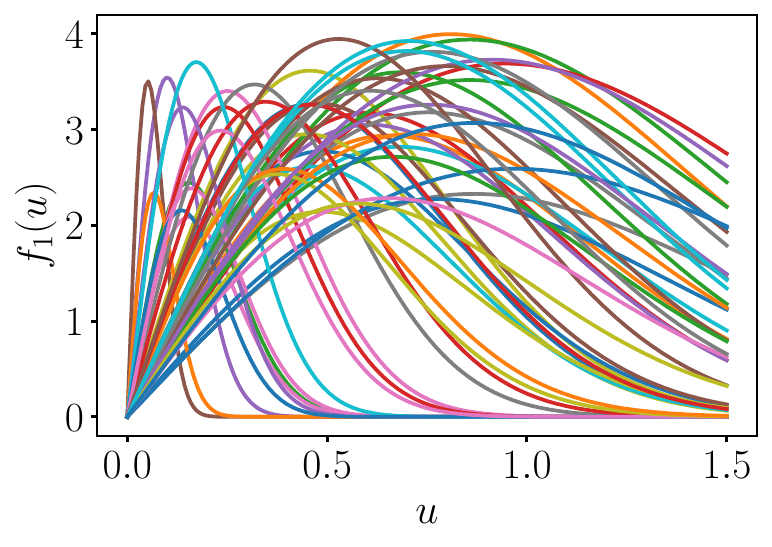}
	\end{subfigure}
	\begin{subfigure}[t]{0.325\linewidth}
		\centering
		\includegraphics[width=\linewidth]{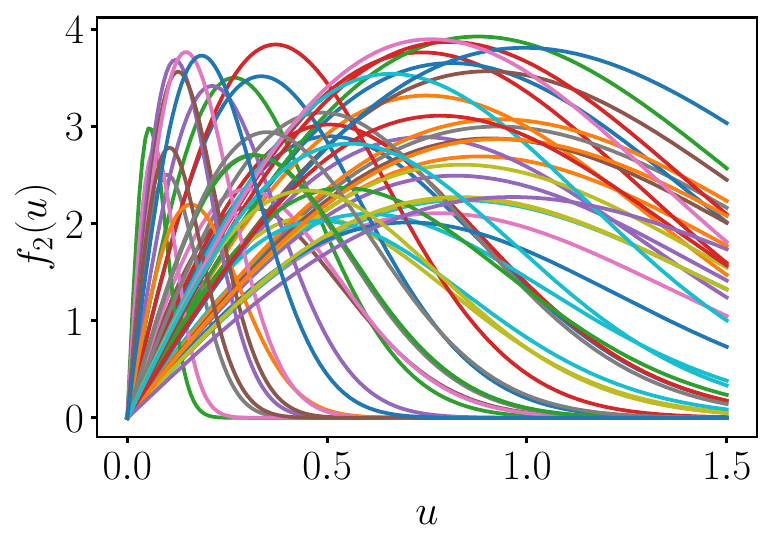}
	\end{subfigure}
	\begin{subfigure}[t]{0.325\linewidth}
		\centering
		\includegraphics[width=\linewidth]{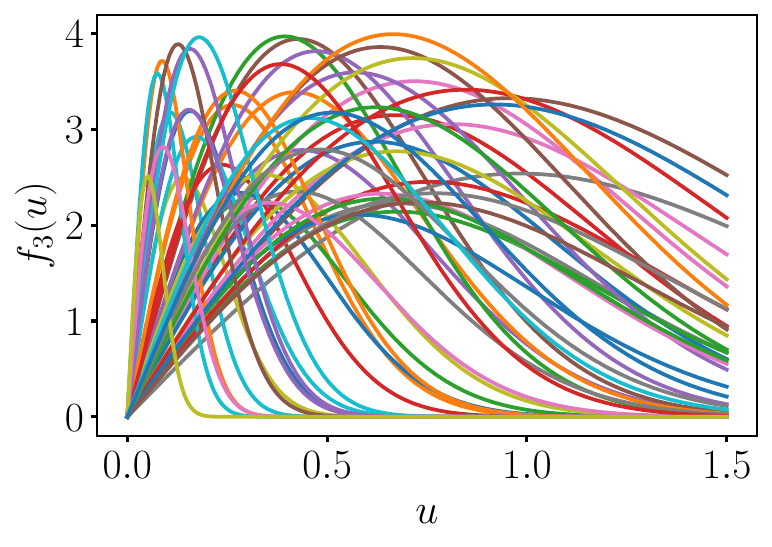}
	\end{subfigure}

	\medskip
	
	\includegraphics[width=0.37\textwidth]{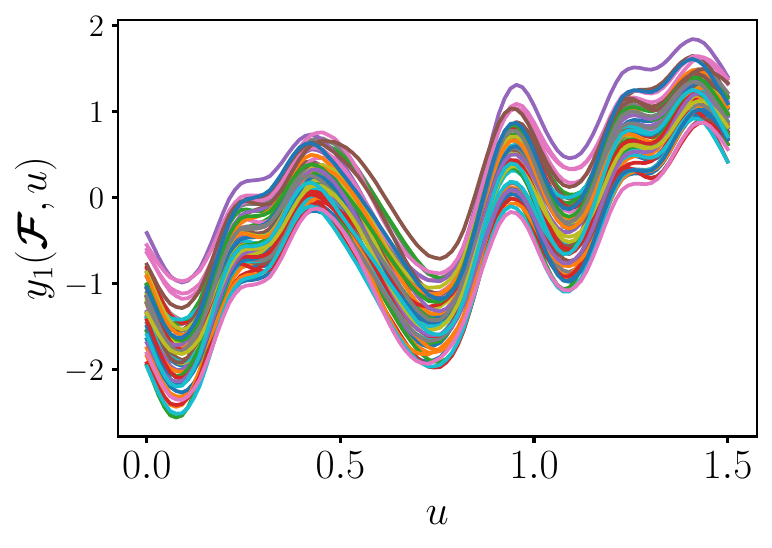}
	\includegraphics[width=0.37\textwidth]{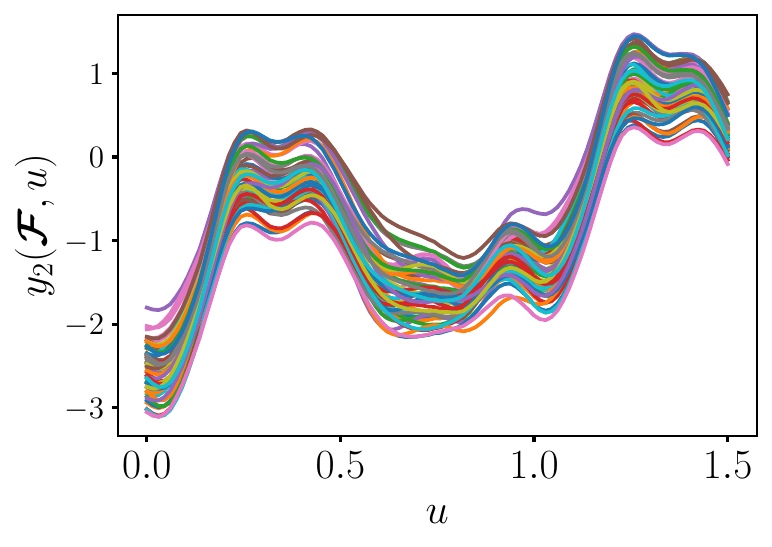}
	
	\caption{Rayleigh-based dataset. The panels show: (top) Samples of the three functional inputs $\bm{\mathcal{F}} = (f_1, f_2, f_3)$, randomly generated as Rayleigh-shaped functions $h_\rho$, and (bottom) the corresponding sampled outputs  $y_1$ and $y_2$ generated by the baseline MTGP.}
	\label{fig:functional_inputs_outputs}
\end{figure}

Figure~\ref{fig:functional_inputs_outputs} shows the corresponding outputs and input profiles generated in this synthetic experiment. For this example, the dataset accounts for  two tasks, i.e., $S=2$. The inter-task correlations are encoded through the  matrix
$
\mathbf{K}_{\mathcal{S}} =
\begin{pmatrix}
	1.00 & 0.85 \\
	0.85 & 1.00
\end{pmatrix}$,
which reflects a strong positive dependence between both tasks. 
The two tasks generated for each set of functional covariates are observed at $n_u = 100$ points, resulting in a total number of
$
n = S \, n_f \, n_u = 10^5$
observations. The two outputs-specific
$y_1(\bm{\mathcal{F}}, u)$ and $y_2(\bm{\mathcal{F}}, u)$
are jointly generated by the baseline MTGP. 

\subsubsection{Cross-validation test}

The predictive performance of the model is assessed using a leave-one-out (LOO)
cross-validation strategy. At each fold, one entire functional replicate is removed
from the training set and used as the test input, while the model is retrained on the
remaining \(n_f - 1\) replicates. To ensure consistency across folds, all functional inputs, whether in the training or test set, are encoded with the same discretization plus PCA procedure as described above. Once the model is  trained, the MTGP predictive mean and variance are evaluated at all \(n_u\) scalar
points of the test curve.
This procedure is repeated until each functional replicate has been used once as the test case, yielding an unbiased estimate of the predictive performance.
Global accuracy is quantified by aggregating the \(Q^2\) scores over all LOO folds.

\begin{figure}[t!]
	\centering
	
	\hspace{3ex}
    \begin{minipage}{0.31\textwidth}
        \centering
 		{\footnotesize Scenario 45}
 	\end{minipage}
 	\begin{minipage}{0.31\textwidth}
        \centering
 		{\footnotesize Scenario 140}
 	\end{minipage}
 	\begin{minipage}{0.31\textwidth}
        \centering
 		{\footnotesize Scenario 429}
 	\end{minipage}
    
    \medskip
	
	\includegraphics[width=\linewidth]{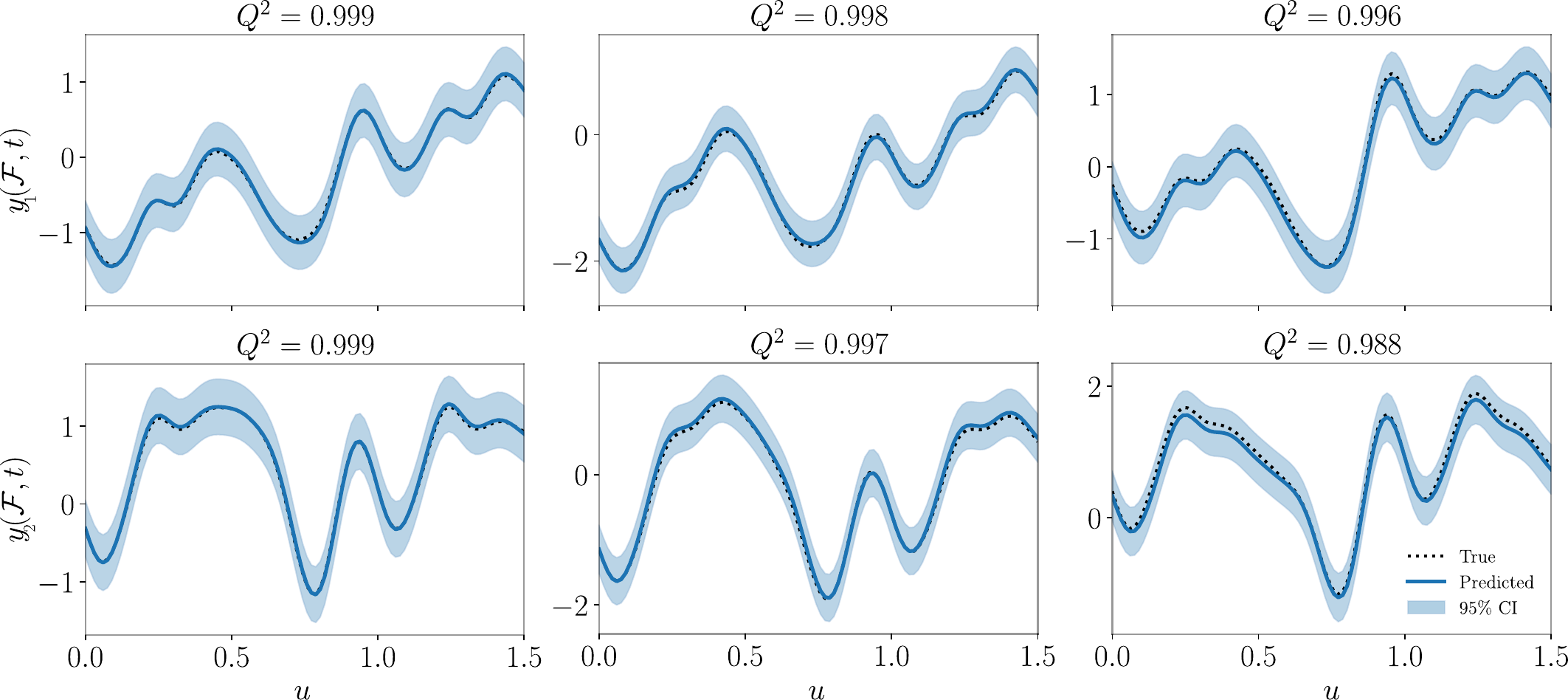}
	
\caption{
    LOO MTGP predictions for the test scenarios 45, 140, and 429 (displayed by columns) of the Rayleigh-based dataset. These scenarios correspond to the ``best'', ``average'', and ``worst'' predictive cases, respectively. The ground truth is shown in black, while predictions are shown in blue with $95\%$ confidence intervals.
}

	\label{fig:predictions_16_45_88}
\end{figure}

Before analyzing the quantitative LOO results, we first present in Figure~\ref{fig:predictions_16_45_88} predictions for three scenarios to enable a visual comparison between the true simulator responses and the MTGP posterior mean, along with the associated \(95\%\) credible intervals. Specifically, we display predictions for Scenarios~45, 140, and~429, representing the ``best'', an ``average'', and the ``worst'' cases in terms of the \(Q^2\) values. We observe that Scenarios~45 and 140 exhibit high predictive accuracy (\(Q^2 \ge 0.997\)) for both tasks, with a correct agreement between predicted and true curves and well-calibrated credible intervals. In contrast, although Scenario~429 yields lower $Q^2$ values (\(Q^2 \ge 0.988\)), both predictions and credible intervals remain accurate and successfully capture the output dynamics.

\begin{figure}[t!]
    \centering
    \begin{minipage}{0.45\textwidth}
        \centering
        \vspace{-2ex}
        \includegraphics[
            width=\linewidth,
            trim=1.5mm 2.5mm 2.5mm 2.5mm,
            clip
            ]{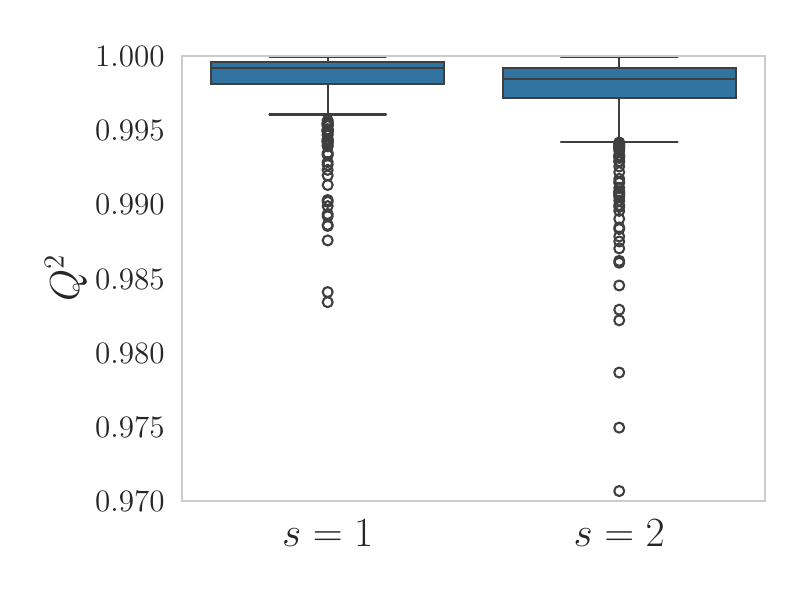}
        \caption{$Q^2$ boxplots computed over the 500 LOO cross-validation replicates used in the Rayleigh-based dataset. Results are shown for each output $y_s$ for $s \in \{1, 2\}$.}
        \label{fig:boxplot_Q2}
    \end{minipage}  
    \hspace{2ex}
    \begin{minipage}{0.43\textwidth}
        \centering
        \includegraphics[
            width=\linewidth,
            trim=2.5mm 2.5mm 2.5mm 2.5mm,
            clip
            ]{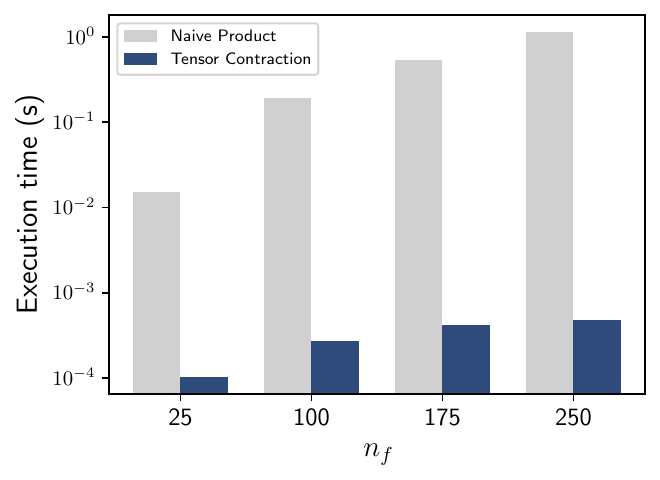}
        \caption{Runtime for solving the linear system $\mathbf{L} \boldsymbol{\alpha} = \boldsymbol{y}$ in the Rayleigh-based dataset. The time is measured for different numbers of functional dimensions $n_f$.}
        \label{fig:boxplot_times}
    \end{minipage}  
\end{figure}

Figure~\ref{fig:boxplot_Q2} reports now the distribution of the 500 LOO \(Q^2\) scores for the two tasks. We can note that the resulting distributions are highly concentrated, with median values close to one, indicating a stable and robust predictive behavior of the MTGP across all test replicates. Regarding the CA results, the model yielded values equal to one in all cases, indicating that the credible intervals always covered the target curves. We have generally noticed that, as observed in Figure~\ref{fig:predictions_16_45_88}, these intervals closely follow the output dynamics while not being overly conservative.

\subsubsection{Computational benchmark}
\label{subsec:runtime-benchmark}

To assess the computational benefit of exploiting Kronecker algebra, we compare two
implementations for applying a Kronecker-structured linear operator to a vector.
More precisely, we consider the computation of $
\boldsymbol{\alpha}$. The first approach, referred to as the naive implementation, explicitly forms the
full Kronecker product
$
\mathbf{L}
=
\mathbf{L}_{\mathcal S}\otimes \mathbf{L}_f\otimes \mathbf{L}_u
\in \mathbb{R}^{n\times n}$,
and computes the matrix--vector product
$\boldsymbol{\alpha}$ 
using standard dense linear algebra routines.
This strategy incurs a quadratic cost in $n$ for each application of the operator,
as well as prohibitive memory requirements.

The second approach corresponds to the proposed tensorized implementation.
In this case, the Kronecker product matrix $\mathbf{L}$ is never formed explicitly.
Instead, the vector $\boldsymbol{y}$ is reshaped into a tensor
$\mathcal{Y}\in\mathbb{R}^{S\times n_f\times n_u}$,
and the action of $\mathbf{L}^{-1}$ is computed via a sequence of mode-wise triangular
solves along the temporal, functional, and task dimensions described in Eq.~\eqref{eq:modewise-prod}.

Both implementations are benchmarked using the same factors
$(\mathbf{L}_{\mathcal S}, \mathbf{L}_f, \mathbf{L}_u)$.
We measure the runtime required to compute
$\boldsymbol{\alpha}$ for increasing functional dimensions
$n_f \in \{25, 100, 175, 250\}$,
averaged over 50 repetitions. Figure~\ref{fig:boxplot_times} compares the execution times of the naive Kronecker
implementation and the tensor contraction approach as a function of the functional
dimension \(n_f\).
While both runtimes increase with \(n_f\), the naive implementation scales much more
steeply, resulting in execution times one to two orders of magnitude larger than those
of the tensorized formulation over the considered range.
This demonstrates the clear computational advantage of tensor contractions over
explicit Kronecker products, independently of covariance estimation or Cholesky
factorization.

\subsection{Application to a riveted mechanical assembly}
\label{sec:mechanical_application}

The objective of this application is to study the dynamic response of a
multi-material riveted assembly subject to uncertainties in the constitutive
behavior of the connectors. Because of the relatively recent adoption of self-piercing riveting, the material
properties remain imperfectly characterized.
Using a limited set of high-fidelity numerical simulations, we assess the ability
of the proposed MTGP to model the relationship between the
force--displacement responses of the connectors and that of the assembly, along with the ability of this model to quantify the  variability in the structural response.

\subsubsection{Mechanical assembly and numerical model}
\begin{figure}[t!]
    \centering
    \includegraphics[width=0.85\textwidth]{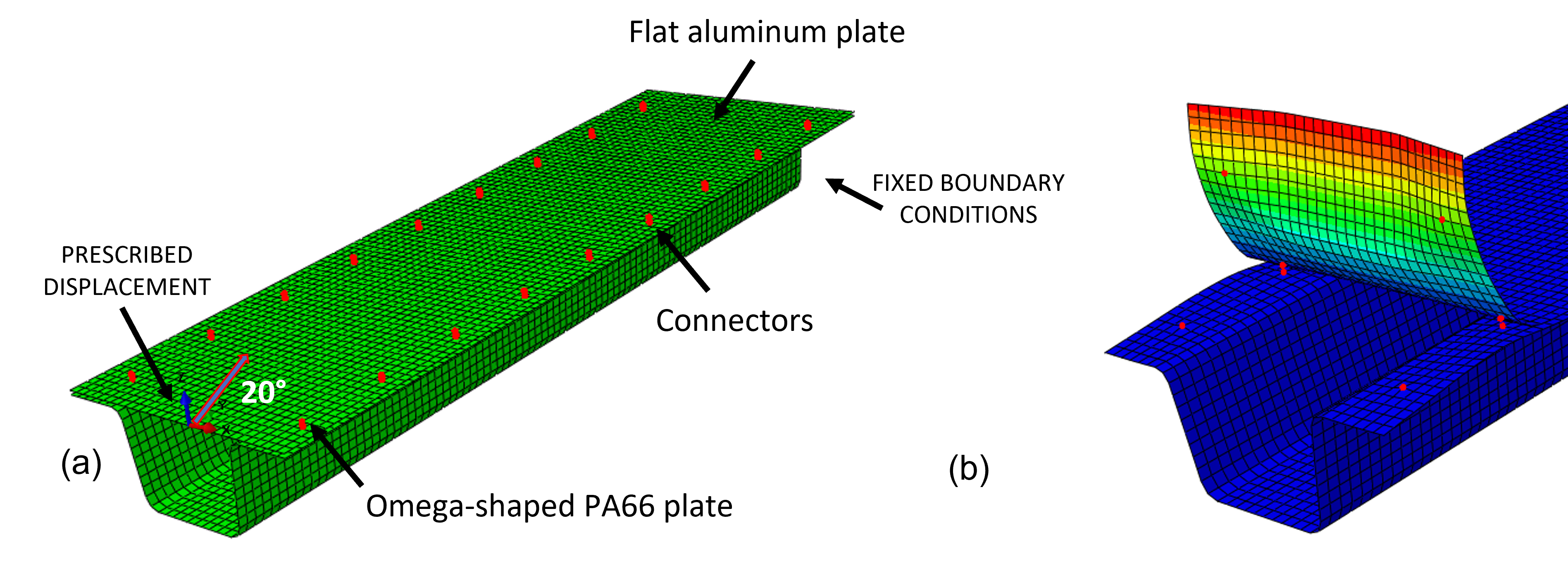}
    \caption{
    Numerical model of the assembly. The panels show: (a) the finite element mesh and simulation conditions, and (b) the deformed shape at the end of the simulation.
    }
    \label{fig:spr_combined}
\end{figure}

The mechanical structure consists of a flat aluminum plate connected to an omega-shaped  plate in polyamid 66 (PA66) by two rows of nine self-piercing rivets. Such assemblies are commonly found in transportation vehicles, especially in the automotive industry.  The physical simulations of the structure are carried out with the finite element model shown in Figure \ref{fig:spr_combined}a. 
All the degrees of freedom of the nodes located on the right side of the PA66 omega-shaped profile are fixed. 
A prescribed displacement is applied to the plate, rotated by $-20^\circ$ around the $x$-axis. The deformed shape of the assembly at the end of the numerical simulation is shown  in Figure \ref{fig:spr_combined}b. Further details on the finite element simulation of the assembly can be found in \citep{Leconte2020}.

\subsubsection{Database description}

\begin{figure}[t!]
	\centering
	\includegraphics[width=0.31\textwidth]{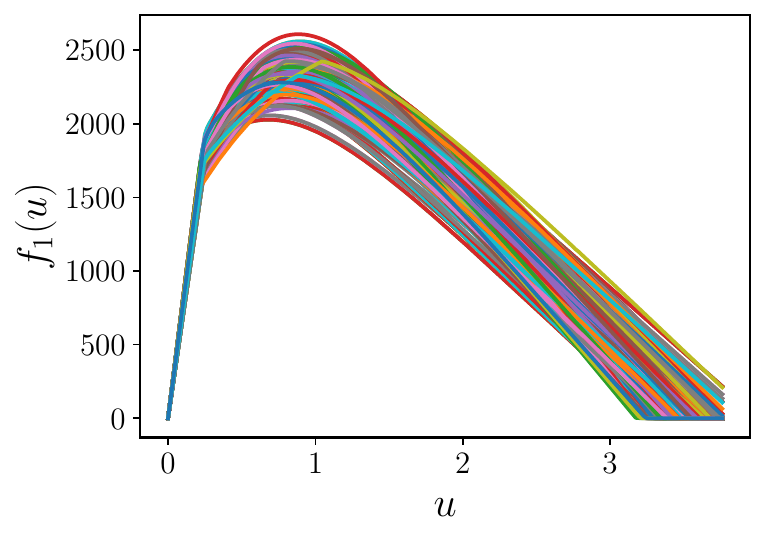}
	\includegraphics[width=0.31\textwidth]{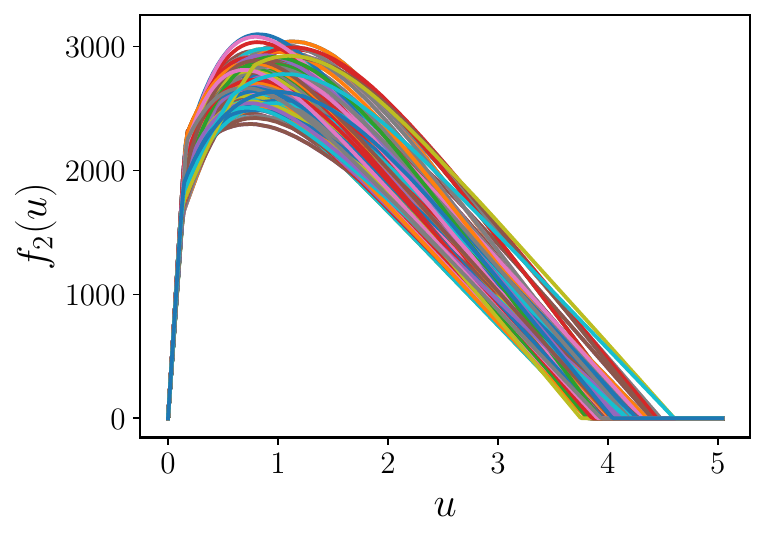}
	\includegraphics[width=0.31\textwidth]{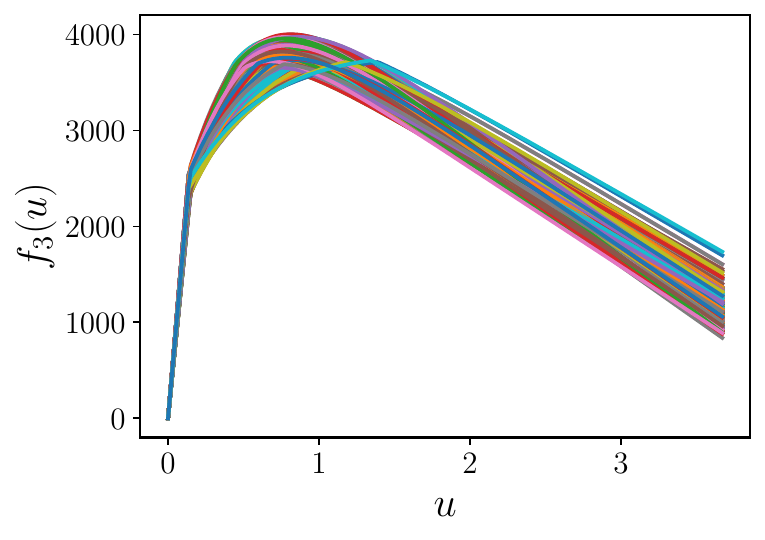}

    \medskip
    
	\includegraphics[width=0.35\textwidth]{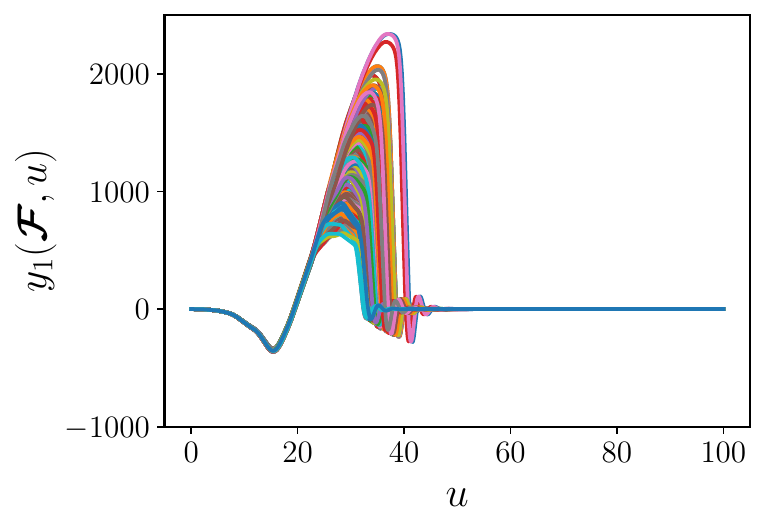}
	\includegraphics[width=0.35\textwidth]{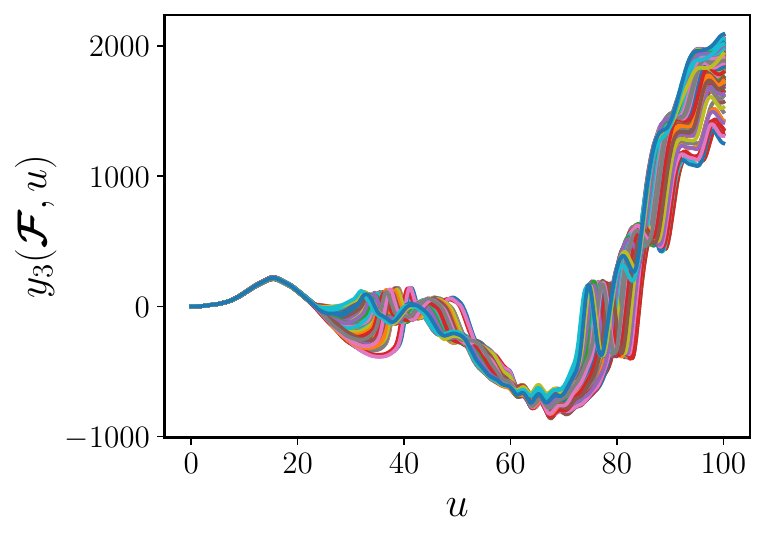}
    \vspace{-1.8em}
    
	\includegraphics[width=0.35\textwidth]{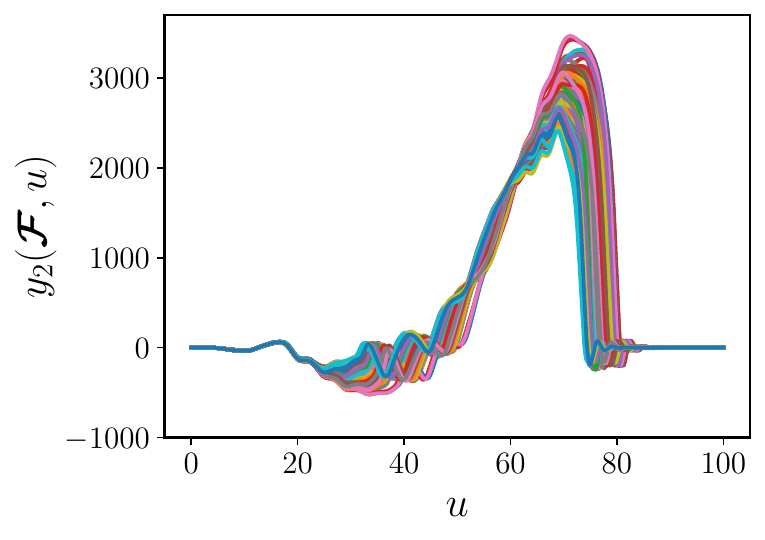}
	\includegraphics[width=0.35\textwidth]{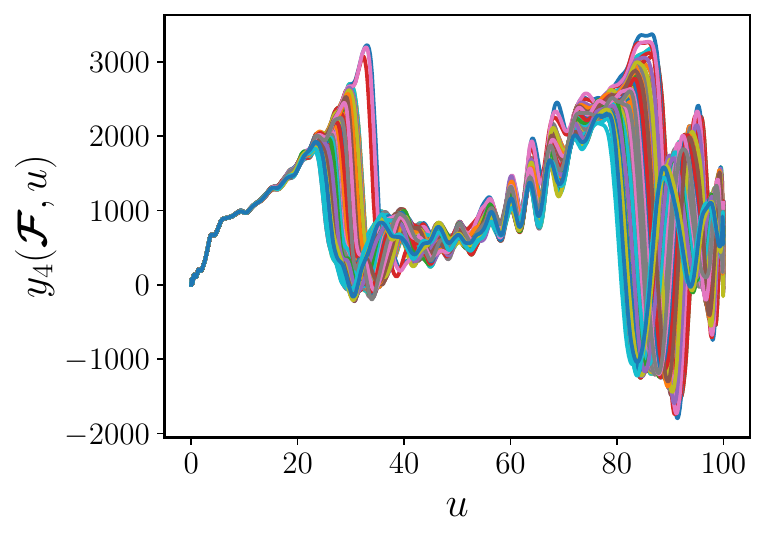}

	\caption{Force-displacement profiles of the riveted assembly. The panels show: (top) the force-displacement functional inputs $\bm{\mathcal{F}} = (f_1, f_2, f_3)$ associated with the connectors, and (bottom) the corresponding outputs $y_s$, for $s \in \{1,\dots,4\}$, associated with the mechanical assembly, each corresponding to a distinct observation location. The horizontal axis represents the prescribed displacement $u$ expressed in millimeters ($mm$), while the vertical axis represents the corresponding force $y_s$ expressed in Newtons ($N$).
}

	\label{fig:functional_inputs_outputs_meca}
\end{figure}

The connectors are mechanically characterized by three functional data $\bm{\mathcal{F}}= (f_1(u), f_2(u), f_3(u))$ (see Figure \ref{fig:functional_inputs_outputs_meca}, top panel), corresponding to the force–displacement responses under pure tension, combined tension–shear, and pure shear loadings, respectively. The different samples reflect material-related uncertainties associated with elastic, plastic and damage behaviors, and constitute the functional covariates used as inputs to the MTGP model. When sampling the inputs, the same triplet 
$\bm{\mathcal{F}}$ is assigned to all connectors of the assembly.  
In other words, every connector within the structure receives the same functional 
inputs $(f_1,f_2,f_3)$ as it is assumed that it is made of the same material. Seventy-eight  samples are employed to explore the design space of these functional inputs.  
For each sampled triplet $\bm{\mathcal{F}}$, a full nonlinear explicit finite 
element simulation is run with the Abaqus software  to extract the corresponding force–displacement responses of the assembly. 
Each finite element simulation requires several hours of computation, which strongly motivates the development of surrogate models capable of delivering predictions of the connector and assembly responses under functional inputs.
Four output functional responses  
$(y_1(\bm{\mathcal{F}},u),\, y_2(\bm{\mathcal{F}},u),\, y_3(\bm{\mathcal{F}},u),\, y_4(\bm{\mathcal{F}},u))$  
(Figure~\ref{fig:functional_inputs_outputs_meca}, bottom panel) are considered, corresponding to the force-displacement responses of the left 
connectors in the first three rows and to the global structural response.  
In a compact form, the resulting dataset is expressed as the functional mapping
\[
    \bm{\mathcal{F}} = (f_1(u), f_2(u), f_3(u))
    \;\mapsto\;
    (y_1(\bm{\mathcal{F}}, u),\, y_2(\bm{\mathcal{F}}, u),\, y_3(\bm{\mathcal{F}}, u),\, y_4(\bm{\mathcal{F}}, u)).
\]

In the MTGP model, the similarity between functional inputs is modeled through the
functional kernel \(k_f\), implemented as a Matérn \(5/2\) covariance acting on the
latent functional representations.
The corresponding length-scales are initialized to \((20,20,20)\), enforcing a smooth
prior across the functional dimensions before being optimized during training.
Displacements are modeled using a scalar kernel \(k_u\), defined as a Matérn \(5/2\)
covariance function.
The initial length-scale is set to \(10\) mm and constrained to the interval
\([1,\,50]\) mm to ensure numerical stability and physically meaningful spatial
variability.
These choices are consistent with the displacement operating range observed in the output profiles shown in
Figure~\ref{fig:functional_inputs_outputs_meca} (bottom panel).

\begin{figure}[t!]
	\centering
	\includegraphics[width=0.49\textwidth]{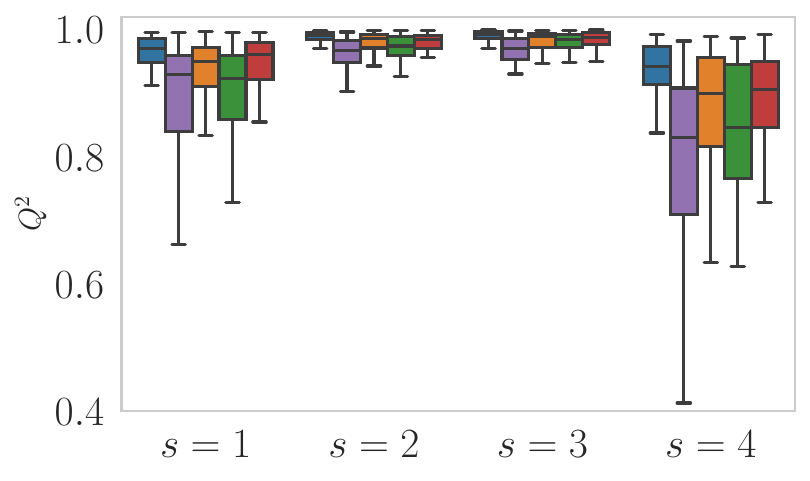}
	\includegraphics[width=0.49\textwidth]{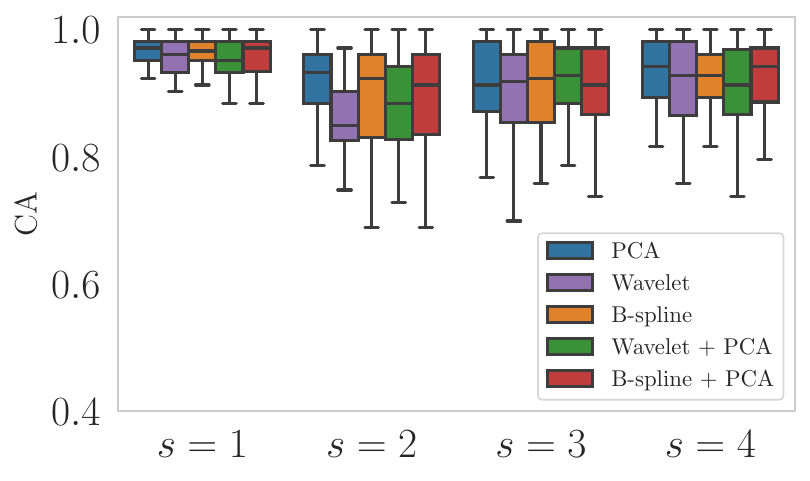}
\caption{
	Predictive accuracy ($Q^2$, left) and uncertainty quantification
	($\operatorname{CA}$ with $\delta = 1.96$, right) obtained with
	$n_f = 78$ training cases and $n_{\text{test}} = 50$ test cases,
	across the four tasks ($s = 1,\dots,4$), using different functional
	dimensionality reduction methods.
}

	\label{fig:gp_perf_tasks_78}
\end{figure}

\begin{figure}[t!]
	\centering
	
	\hspace{5ex}
    \begin{minipage}{0.305\textwidth}
        \centering
 		{\footnotesize Scenario 23}
 	\end{minipage}
 	\begin{minipage}{0.305\textwidth}
        \centering
 		{\footnotesize Scenario 40}
 	\end{minipage}
 	\begin{minipage}{0.305\textwidth}
        \centering
 		{\footnotesize Scenario 29}
 	\end{minipage}
 	
 	\medskip
	
	\includegraphics[width=\textwidth]{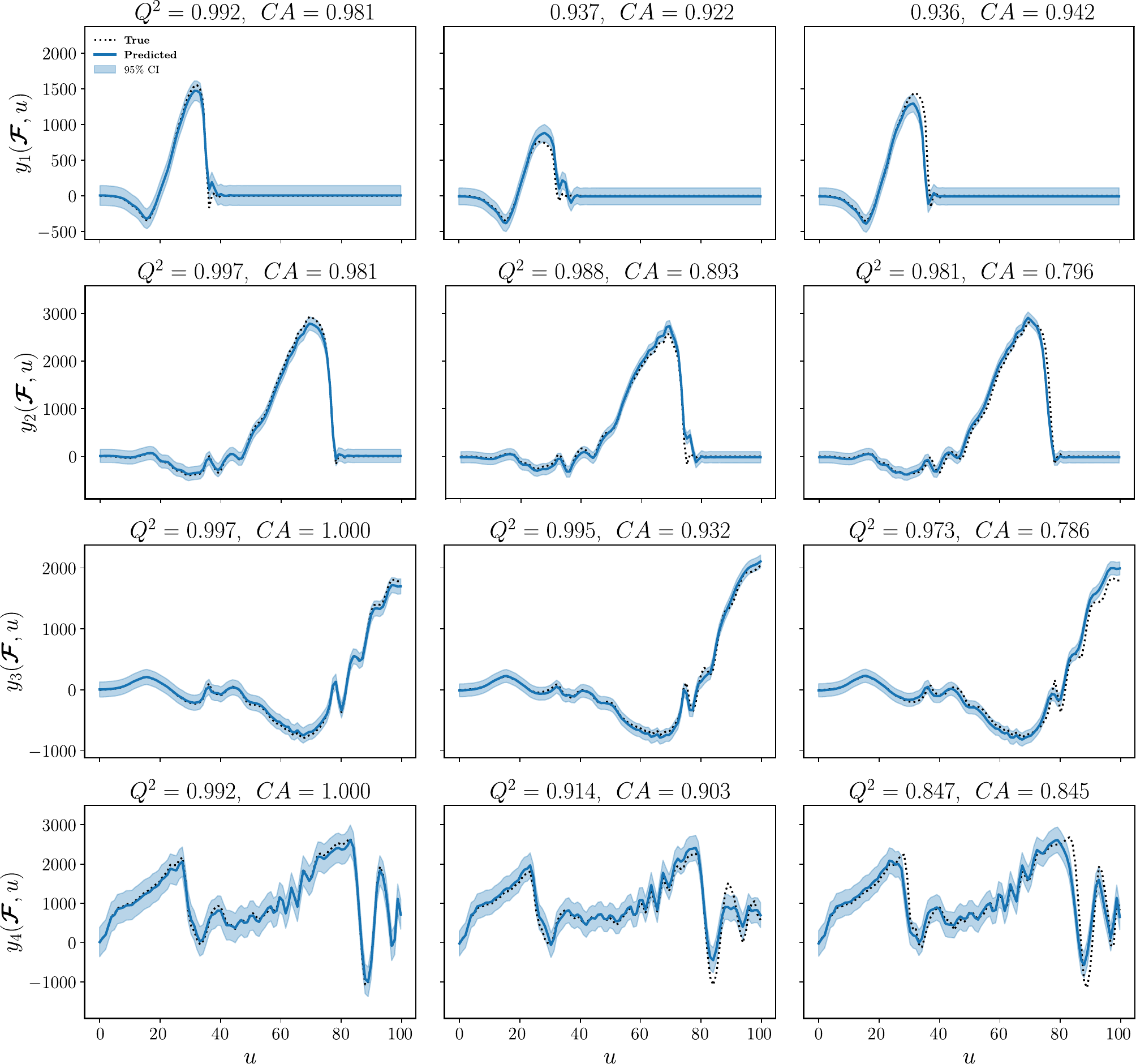}

    \caption{MTGP predictions for scenarios test ~23,~40, and~29 (displayed by columns) of the riveted  assembly. These scenarios represent the best, the average, and the worst predictive cases, respectively, in terms of the $Q^2$ criterion.
   The ground truth is the black dotted line, while the predictions are plotted in blue with the 95\% confidence intervals
    .}
    
	\label{fig:gp_pred_examples}
	
	\bigskip
	
	\includegraphics[width=0.49\textwidth]{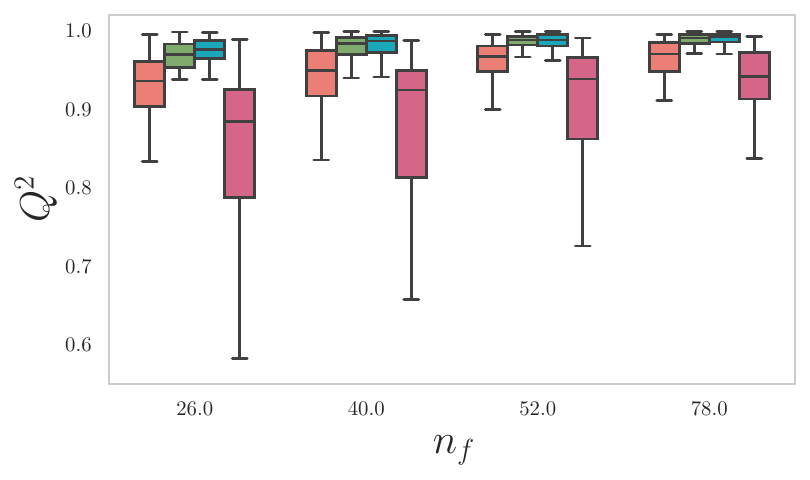}
	\includegraphics[width=0.49\textwidth]{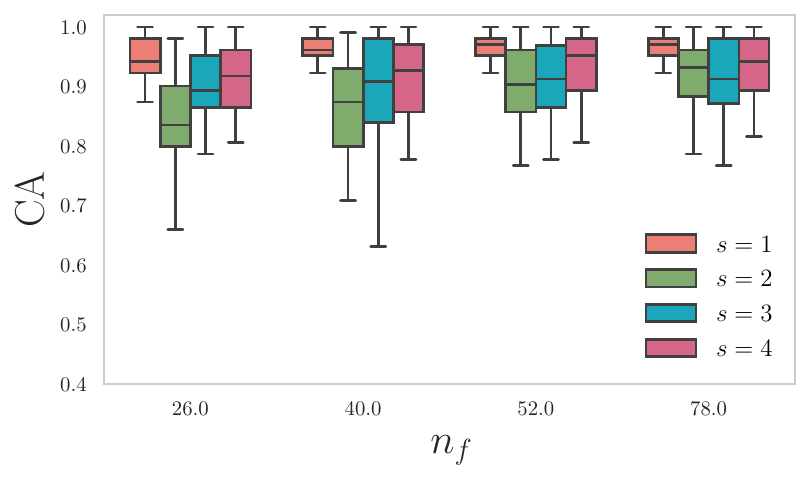}
	
\caption{
	Influence of the training set size on predictive accuracy ($Q^2$, left) and uncertainty
	calibration (CA, right) across the four tasks ($s=1,\dots,4$)  for the riveted assembly test case.
}
	\label{fig:gp_perf_training_size_extra}
\end{figure}

\subsection{Results and discussion}
\label{subsec:gp_perf_training_size}

\subsubsection{Impact of functional dimensionality reduction on predictive performance}
The purpose of this experiment is to assess the influence of the functional dimensionality reduction strategy on the predictive performance of the MTGP model. Figure~\ref{fig:gp_perf_tasks_78} reports results obtained with
$n_f = 78$ training functions and $n_{\text{test}} = 50$ test functions,
for the four tasks ($s \in \{1,2,3,4\}$),
comparing PCA, wavelet, B-spline, and hybrid wavelet+PCA and B-spline+PCA encodings.  From a predictive accuracy perspective, it can be observed that the direct
PCA-based encoding  yields the highest performance across all tasks.
This trend is particularly visible for the most challenging outputs, where PCA
exhibits higher median accuracy together with reduced dispersion, suggesting a more stable predictive behavior.
Moreover, the PCA predictive intervals, as measured by the CA, remain close to the nominal coverage level (about 0.95 for $\delta=1.96$),  indicating an uncertainty quantification that matches well the model assumptions.
Among the projection-based approaches, the B-spline encodings achieve results that
are close to those obtained by the PCA. With the B-splice, the performance remains relatively stable across tasks, both in terms of
predictive  and uncertainty accuracy, which suggests that B-splines
provide an effective low-dimensional representation of the functional inputs.
By contrast, the wavelet-based representation leads to lower predictive accuracy
and increased variability, especially for outputs $y_1$ and $y_4$, which are
characterized by more complex response patterns.
Applying PCA to the wavelet coefficients partially mitigates these limitations:
the Wavelet+PCA strategy improves both the median accuracy and the stability of the
predictions compared to raw wavelet features.
However, its performance remains below that of direct PCA and
B-spline encodings.
Similarly, the B-spline+PCA variant induces only limited changes relative to the
original B-spline representation, indicating that the B-spline coefficients are
already weakly correlated and sufficiently structured, so that additional linear
compression provides little improvement.
Overall, these results indicate that direct PCA-based functional encoding offers
a favorable compromise between predictive accuracy, robustness, and uncertainty
quantification across the four tasks.
For this reason, PCA is retained as the functional dimensionality reduction method
in all subsequent experiments.

To further illustrate these results, we now examine the MTGP predictions over test scenarios. In Figure~\ref{fig:gp_pred_examples}, prediction
examples are reported for three selected test scenarios (23, 40, and 29) which
correspond, respectively, to the ``best'', ``average'', and ``worst'' predictive cases according
to the \(Q^2\) criterion.
The results are shown for all four tasks. 
From the figure, one can observe that the first three responses
(\(y_s\), $s = 1,2,3$) are generally reconstructed with good accuracy, whereas the
fourth task (\(s = 4\)) displays a comparatively lower level of agreement.
This behavior is consistent with the trends previously seen in
Figure~\ref{fig:gp_perf_tasks_78}.
 This difference is not unexpected, since the output \(y_4\)
corresponds to a more global structural response that aggregates complex local
effects and involves stronger nonlinearities.
Despite this increased level of difficulty, the MTGP remains able to capture the
main temporal evolution of the global $y_4$ in a generalization regime. 

For completeness, additional prediction results obtained with alternative
functional projection methods are provided in~\ref{app:comp_results}.

\subsubsection{Influence of the training set size on prediction and uncertainty calibration}
We now investigate how the size of the training set influences both the predictive accuracy and the uncertainty quantification of the MTGP model. To this end, the number of training functions is gradually increased ($n_f \in \{26, 40, 52, 78\}$), with each dataset constructed in a nested manner within the larger ones. The corresponding results are reported in Figure~\ref{fig:gp_perf_training_size_extra}. 
As could be expected, increasing the training set size leads to a improvement in predictive performance across all tasks. This evolution is reflected by higher $Q^2$ values together with coverage accuracy progressively approaching the reference level, indicating a gradual improvement in uncertainty calibration as more training information becomes available. 
A closer inspection reveals that, for tasks $y_2$ and $y_3$, a high predictive accuracy is already achieved with relatively small training sets, with $Q^2$ values close to one as soon as $n_f$ exceeds 40. 
However, for these tasks, the corresponding coverage accuracy remain below the levels obtained for $n_f = 78$. Although the predictive mean is accurately captured at an early stage, the associated predictive uncertainty tends to be underestimated, leading to overly narrow prediction intervals. Such underestimation of the uncertainty is related to the estimated inter-task correlation structure, as illustrated by the task correlation matrix $\mathbf{K}_{\mathcal S}$ (Figure~\ref{fig:task_correlation_matrix}), which shows strong dependencies between tasks $y_2$ and $y_3$.
There, information sharing across tasks supports accurate mean prediction even with limited data, while the estimation of task-specific uncertainty remains more sensitive to the size of the training set. As the number of training functions increases, this effect progressively diminishes, resulting in improved uncertainty calibration across all tasks. 

\subsubsection{Assessment of predicted calibration envelopes}
Beyond pointwise accuracy, a reliable surrogate must provide uncertainty estimates
that are consistent with the variability of the true system.
To assess this property, we compare the MTGP predictive behavior with the empirical
dispersion of the test responses through a calibration envelope analysis. For each task~$s$, the collection of true simulator outputs
$\{y_s^{\mathrm{true}}(\bm{\mathcal{F}}_{\star,i},u_{\star})\}_{i=1}^{n_{\mathrm{test}}}$
associated with all test functional inputs is used to construct an empirical
$95\%$ envelope defined by
\begin{equation*}
\operatorname{CI}^{\mathrm{true}}_s
=
\left[
q_{2.5\%}\!\left(
y_s^{\mathrm{true}}(\bm{\mathcal{F}}_{\star,i},u_\star)
\right),
\;
q_{97.5\%}\!\left(
y_s^{\mathrm{true}}(\bm{\mathcal{F}}_{\star,i},u_\star)
\right)
\right],
\label{eq:ci_true}
\end{equation*}
where $q_{2.5\%}$ and $q_{97.5\%}$ denote the empirical $2.5$th and $97.5$th percentiles, respectively.
This envelope represents the  variability of the mechanical response induced by the diversity of material conditions, and serves as a reference for
assessing the consistency of the surrogate predictions. 

At the same evaluation points, the MTGP predictive distribution is characterized
by its conditional posterior mean and variance,
denoted by $m_s(\bm{\mathcal{F}}_{\star,i},u_\star)$ and
$v_s(\bm{\mathcal{F}}_{\star,i},u_\star)$, respectively, as defined in
Eq.~\eqref{eq:post-mean-final} and Eq.~\eqref{eq:post-var-final}.
Two complementary aspects of uncertainty are then examined. First, we quantify the variability of the predicted response across test functional inputs
through the dispersion of the MTGP predictive means
$\{m_s(\bm{\mathcal{F}}_{\star,i},u_\star)\}_{i=1}^{n_{\mathrm{test}}}$.
This variability is summarized by the empirical $95\%$ envelope
\begin{equation*}
\operatorname{CI}_{m_s}
=
\left[
q_{2.5\%}\!\left(m_s(\bm{\mathcal{F}}_{\star,i},u_\star)\right),
\;
q_{97.5\%}\!\left(m_s(\bm{\mathcal{F}}_{\star,i},u_\star)\right)
\right].
\label{eq:ci_mean}
\end{equation*}
Second, we assess the local epistemic uncertainty conveyed by the predictive variance
$v_s(\bm{\mathcal{F}}_{\star,i},u_\star)$.
For each test input, a pointwise uncertainty bound is defined as
$m_s(\bm{\mathcal{F}}_{\star,i},u_\star)\pm 1.96\sqrt{v_s(\bm{\mathcal{F}}_{\star,i},u_\star)}$.
Aggregating these bounds over all test inputs yields the envelope
\begin{equation*}
\operatorname{CI}_{v_s}
=
\left[
q_{2.5\%}\!\left(
m_s(\bm{\mathcal{F}}_{\star,i},u_\star)
- 1.96 \sqrt{v_s(\bm{\mathcal{F}}_{\star,i},u_\star)}
\right),
\;
q_{97.5\%}\!\left(
m_s(\bm{\mathcal{F}}_{\star,i},u_\star)
+ 1.96 \sqrt{v_s(\bm{\mathcal{F}}_{\star,i},u_\star)}
\right)
\right].
\label{eq:ub_post}
\end{equation*}

\begin{figure}[t!]
    \centering
    \includegraphics[width=0.45\textwidth]{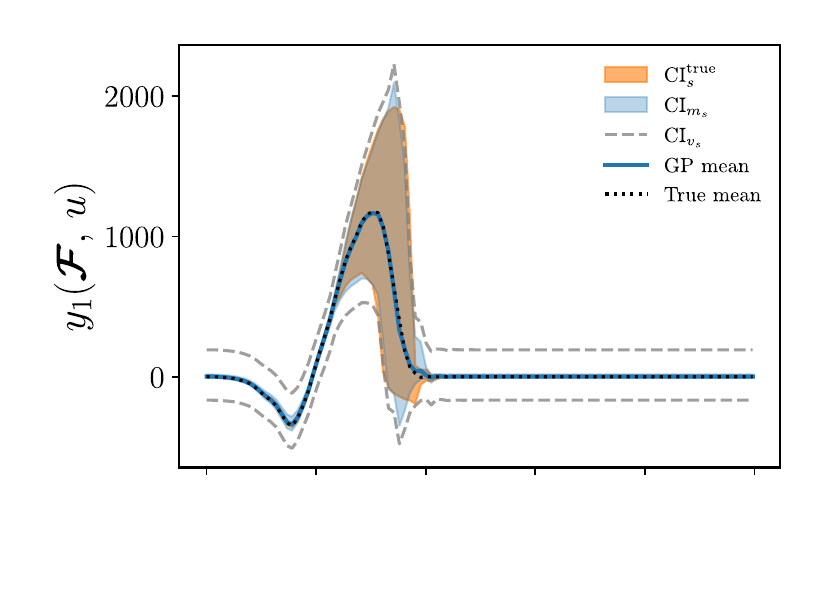}%
    \includegraphics[width=0.45\textwidth]{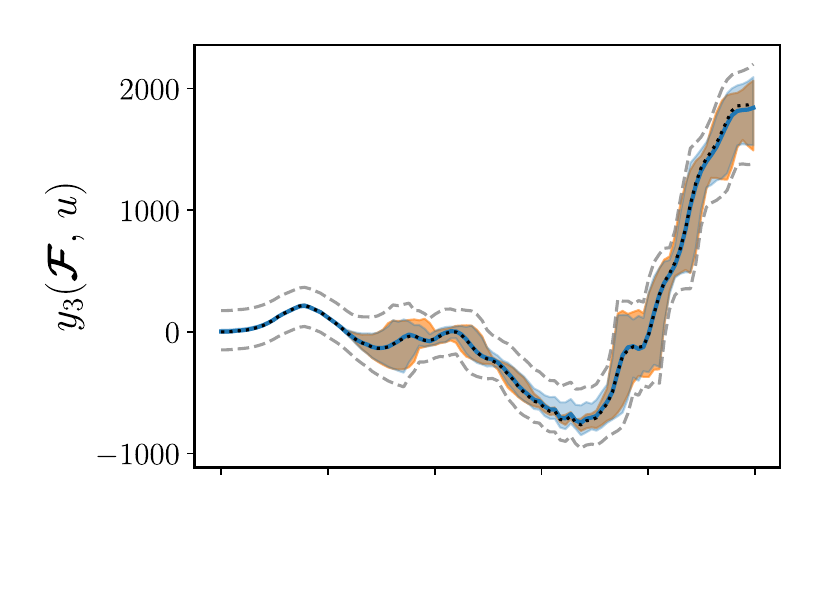}

    \vspace{-2.0em}

    \includegraphics[width=0.45\textwidth]{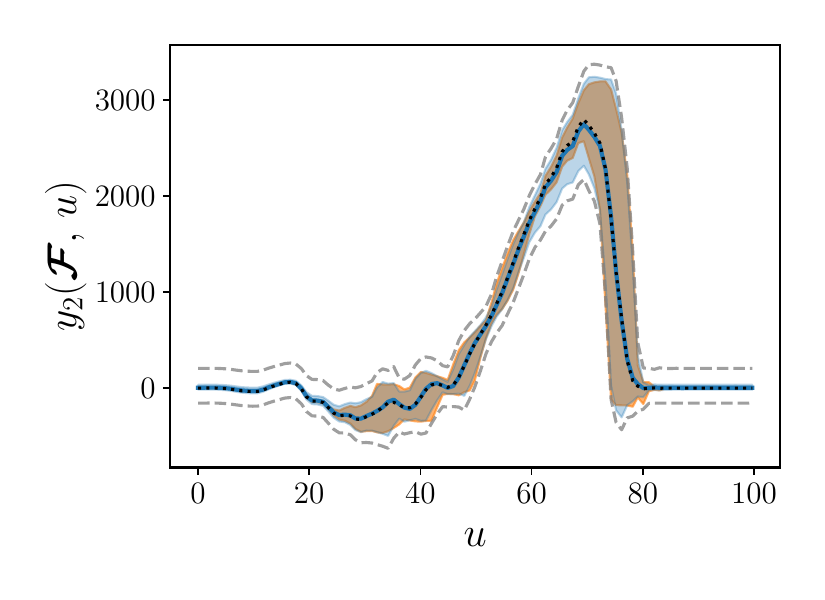}%
    \includegraphics[width=0.45\textwidth]{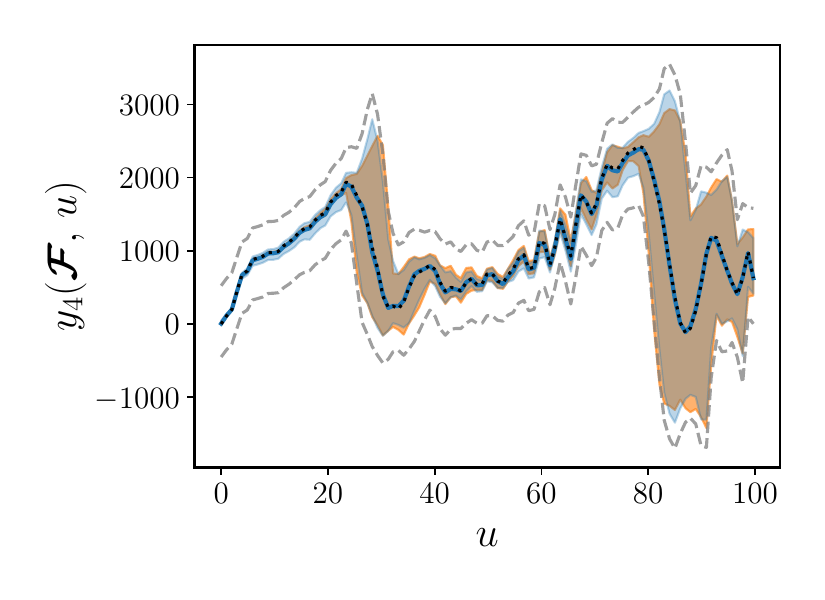}

    \caption{Calibration envelopes for the four MTGP tasks. The orange band shows the empirical $95\%$ envelope of the true simulator responses $\operatorname{CI}^{\mathrm{true}}_s$, the translucent blue region is the predictive envelope from the dispersion of GP mean predictions, $\operatorname{CI}_{m_s}$ and the dashed grey curves denote the $95\%$ posterior credible interval, $\operatorname{CI}_{v_s}$. The solid blue line is the MTGP predictive mean, and the black dotted line represents the empirical mean of the true responses.}
    \label{fig:calibration_envelope}
\end{figure}

Overlaying the three envelopes, $\operatorname{CI}^{\mathrm{true}}_s$, $\operatorname{CI}_{m_s}$ and $\operatorname{CI}_{v_s}$,
provides a visual diagnostic of uncertainty calibration where $\operatorname{CI}^{\mathrm{true}}_s$ serves as the reference for the intrinsic
variability of the mechanical response over the set of test functional inputs. This is done in Figure~\ref{fig:calibration_envelope}. $\operatorname{CI}^{\mathrm{true}}_s$ highlights regions of increased dispersion, that occur in particular in transient
regimes and high-amplitude response zones. The MTGP predictive mean (solid blue curve)  closely follows the empirical mean of the true simulator responses
(black dotted curve), including in regions characterized by rapid temporal variations.
This agreement suggests that the dominant temporal structure of each output is well
captured by the surrogate model.

 The dispersion of the MTGP predictive means across test
functional inputs, summarized by the predictive envelope $\operatorname{CI}_{m_s}$
(blue shaded area) is similar in most regions to that of the empirical envelope, indicating
that the MTGP reproduces, to a large extent, the variability induced by changes in the
functional inputs. Still minor differences between $\operatorname{CI}_{m_s}$ and $\operatorname{CI}^{\mathrm{true}}_s$ can be found at sharp changes in trend. On the other hand, the uncertainty bands derived from the predictive variance,
$\operatorname{CI}_{v_s}$ (dashed grey curves), provide insight into the epistemic
uncertainty of the surrogate model. These bands are wider than the empirical envelope, such conservatism being more pronounced, relatively, when responses are near constant.  Generally however, $\operatorname{CI}_{v_s}$ remains compatible with the scale of the empirical envelope. The larger width of $\operatorname{CI}_{v_s}$ with respect to $\operatorname{CI}^{\mathrm{true}}_s$ and $\operatorname{CI}_{m_s}$ is because the epistemic uncertainties on the inputs are taken into account twice, once through the variance $v_s$, and another time through the quantiles on the inputs.

Overall, these experiments show that the MTGP performs well both in terms of 
predictive accuracy and uncertainty quantification, thereby supporting its use as a
reliable surrogate model for the considered mechanical application.

\subsubsection{MTGP versus single-task GPs}
\label{subsec:MTGP_vs_indGP}

In order to quantify the benefit of explicitly modeling inter-task dependencies,
this section compares the proposed MTGP with functional inputs to its single-task
GP counterpart, which assumes statistical independence among the tasks.
In the single-task setting, each task $s \in \{1,\dots,S\}$ is modeled independently
according to
$
Y_s(\bm{\mathcal{F}},u) \sim \mathcal{GP}\!\big(0,\;
k_{f}(\bm{\mathcal{F}},\bm{\mathcal{F}}')\,k_{u}(u,u')\big)$, where $k_{f}$ and $k_{u}$ capture correlations over the functional and temporal
domains, respectively.
The associated hyperparameter vector is
$
\boldsymbol{\theta}
= \left( \sigma^2,\, \boldsymbol{\ell}_f,\, \ell_{u} \right)$, with $\sigma^2$ denoting the global variance parameter, and
$\boldsymbol{\ell}_f$ and $\ell_{u}$ the length-scale parameters of $k_{f}$ and
$k_{u}$.
The MTGP extends this formulation by introducing an explicit task-covariance matrix
$\mathbf{K}_{\mathcal{S}}$, endowed with its own parameters, so that the only
structural difference between the two models lies in the explicit modeling of
cross-task correlations.

For this comparison, both the single-task GP and the MTGP models are trained under identical
conditions, using the same multi-start optimization strategy and initialization
settings, and with a fixed training size of $n_f = 78$, ensuring a fair and
consistent evaluation. We then compare the predictive accuracy and the uncertainty calibration obtained
with the two approaches.

\begin{figure}[t!]
	\centering
	\includegraphics[width=0.49\textwidth]{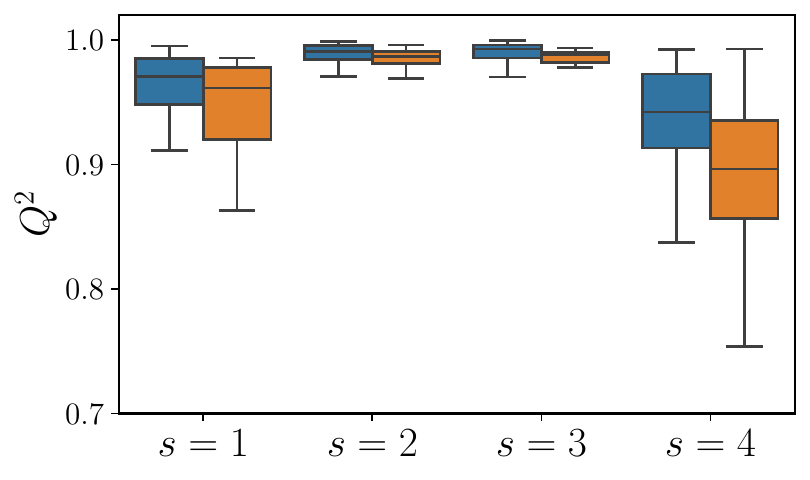}
	\includegraphics[width=0.49\textwidth]{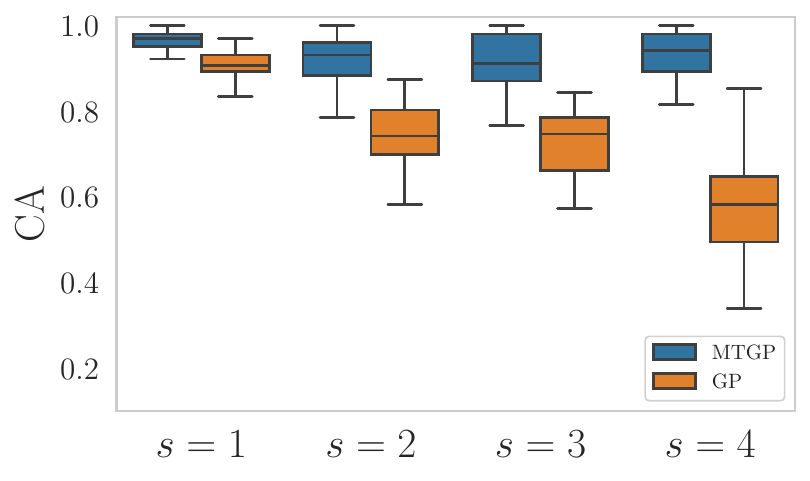}
	\caption{
		Comparison of the prediction  accuracy and uncertainty calibration between MTGP and single-task GPs for a training size of $n_f = 78$. The left panel reports the $Q^2$ distributions across the four output tasks, while the right panel shows the $CA$ coverage accuracy.
	}
	\label{fig:MTGP_vs_indGP_perf}
\end{figure}

As shown in Figure~\ref{fig:MTGP_vs_indGP_perf}, the MTGP consistently
achieves higher $Q^2$ values than the independent GP across tasks,
while maintaining predictive intervals well calibrated around the
nominal $95\%$ level. These results indicate that explicitly modeling cross-task covariance
improves predictive accuracy without degrading uncertainty
quantification. By drawing upon statistical links across outputs, the MTGP provides a more coherent and efficient representation
of the multitask system.

\begin{figure}[t!]
    \centering
    \includegraphics[width=\linewidth]{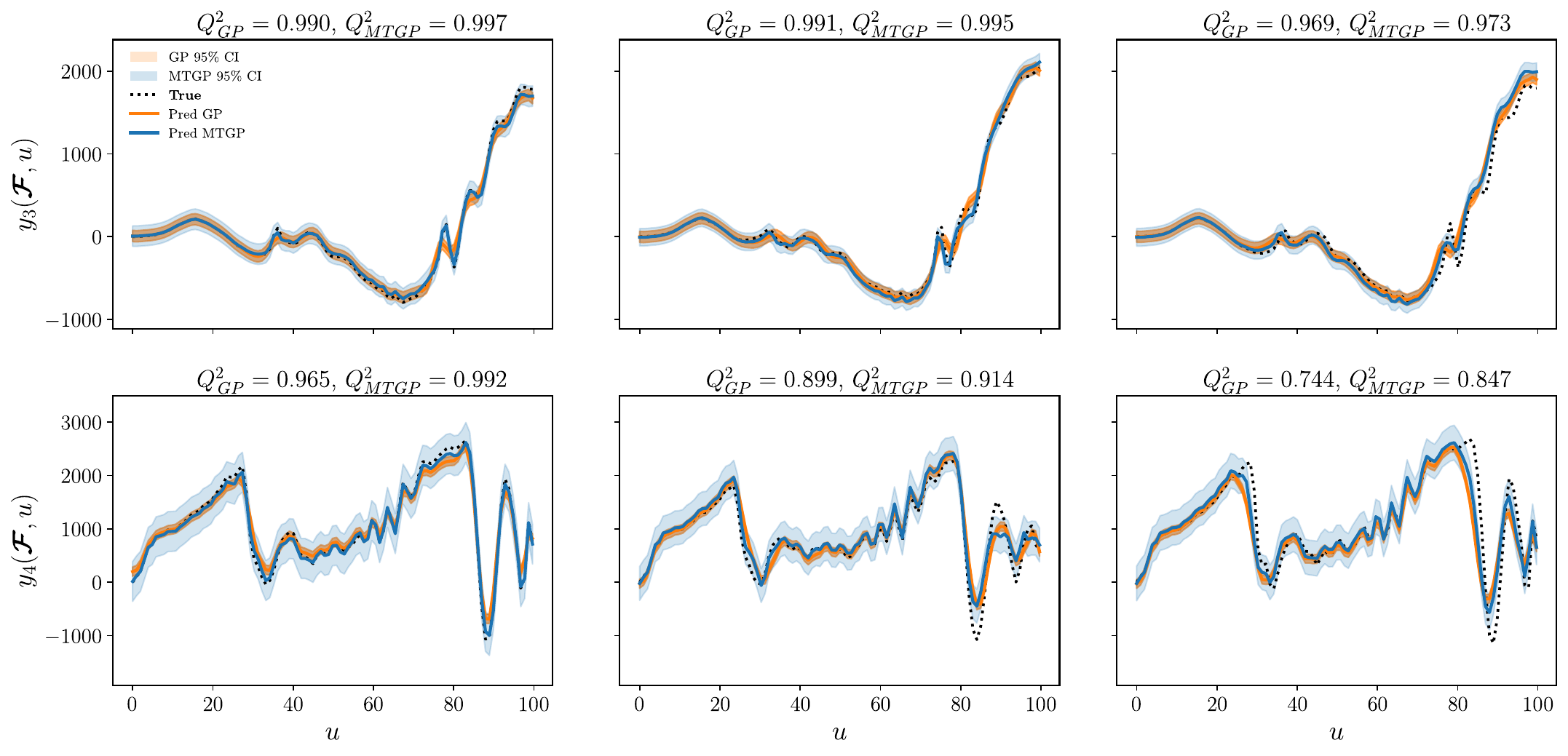}
    \caption{Comparison between MTGP and single-task GP predictions for scenarios 23 (best), 40 (average), and 29 (worst). Results are shown for the tasks (top) $s=3$ and (bottom) $s=4$, corresponding respectively to the tasks for which the MTGP yields the highest and lowest \(Q^2\) values.}
    \label{fig:gp_pred_examplescomparison}
\end{figure}

The example predictions reported in Figure~\ref{fig:gp_pred_examplescomparison} illustrate the behavior of the GP and the MTGP for outputs $s=3$ and $s=4$ under the scenarios~23 (best case), 40 (average case), and 29 (worst case). The MTGP model provides a more robust reconstruction of the temporal dynamics in all the cases. 
Its predictive mean tracks the true trajectories more accurately, including in regions exhibiting rapid variations or  oscillations (most notably for the scenario~29 for $s=4$), whereas the GP shows visible deviations around peaks and regime changes. 
A key difference lies in the uncertainty quantification. The MTGP credible intervals (in blue) are wider than those of the GP, yet they remain regular and consistent with the intrinsic variability of the simulator. 
This controlled width reflects the fact that the MTGP accounts for both local uncertainty and inter–task correlations, thereby preventing the artificial underestimation of predictive variance. 
By contrast, the  GP exhibits a more unstable behavior. 
In some regions, its credible intervals are extremely narrow or even barely visible, indicating a substantial underestimation of uncertainty. Overall, they remain entirely contained within the MTGP intervals, revealing the absence of any mechanism to propagate information across correlated outputs. Such overly optimistic intervals do not reflect the true variability of the responses and lead to unreliable uncertainty estimates. Overall,  these results demonstrate that explicitly modelling inter–task correlations enhances both the fidelity of the predictions and the credibility of the associated uncertainty, whereas independent modelling fails to guarantee these properties.

\begin{figure}[t!]
	\centering
	\includegraphics[width=0.45\textwidth]{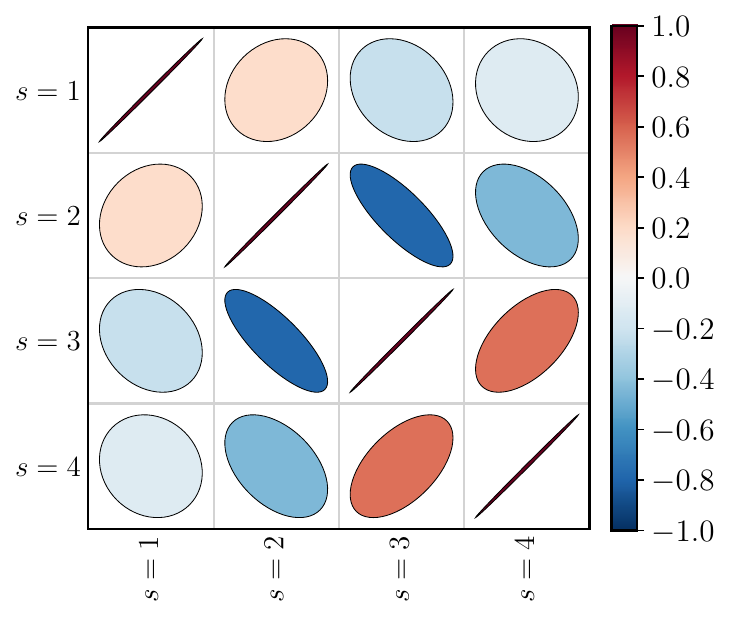}
	\caption{
		Estimated inter-task correlation matrix $\mathbf{K}_{\mathcal{S}}$ for the four tasks involved in the riveted assembly application. Red (respectively blue) ellipses indicate positive (resp. negative) correlations,
with orientation encoding the sign and eccentricity reflecting the magnitude.}
	\label{fig:task_correlation_matrix}
\end{figure}

To quantify the dependencies seen in the data, Figure~\ref{fig:task_correlation_matrix} presents the  task-correlation matrix $\mathbf{K}_{\mathcal{S}}$ estimated after training the MTGP. Substantial correlations are observed most notably a pronounced negative correlation between tasks $s=2$ and $s=3$. This confirms the intuition that the mechanical responses at neighboring rivets are more correlated, thus justifying the use of a multitask formulation to jointly learn the rivets force-displacement relationship.

\begin{figure}[t!]
	\centering
	\includegraphics[width=0.49\textwidth]{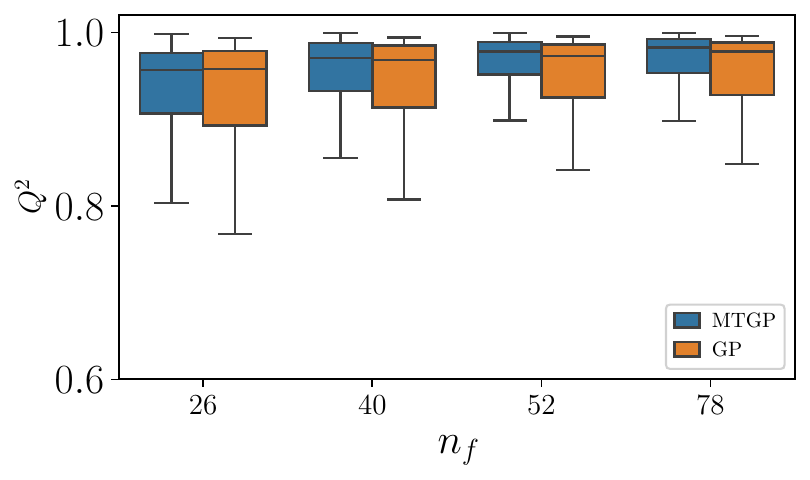}
	\includegraphics[width=0.49\textwidth]{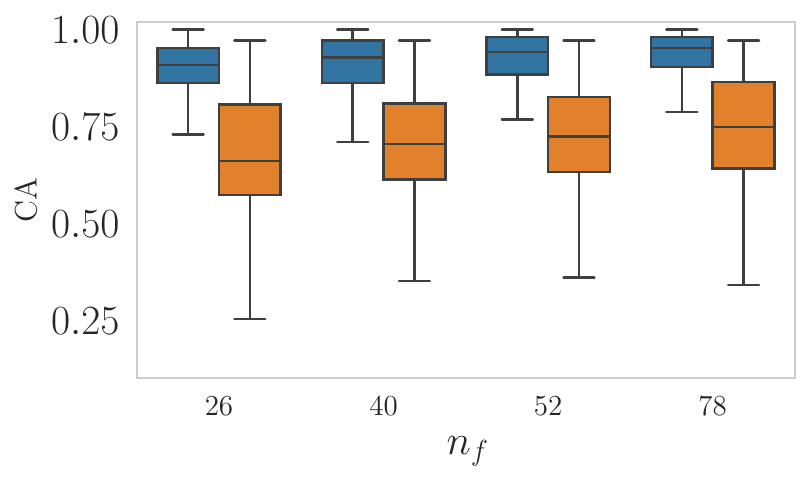}
	\caption{
		Evolution of predictive accuracy and uncertainty calibration with increasing training size 
		($n_f \in \{26, 40, 52, 78\}$) for MTGP and GP models.
	}
	\label{fig:MTGP_vs_indGP_evolution}
\end{figure}

To further compare the dependencies of both single and multi-task GPs on data availability, we analyze how their predictive performances evolve as the number of training functions increases ($n_f \in \{26, 40, 52, 78\}$). 
Figure~\ref{fig:MTGP_vs_indGP_evolution} summarizes this evolution for both $Q^2$ and $\operatorname{CA}$ metrics. 
Each boxplot represents the distribution of these metrics averaged over the four tasks for both models. 
The results show that, as the training set grows, both models improve in terms of accuracy and uncertainty calibration. 
However, the MTGP, which is marginally better than the GP in terms of $Q^2$ when $n_f=26$, leads to faster improvements and rapidly achieves higher $Q^2$ scores and better-calibrated coverage rates, demonstrating its ability to exploit shared information across tasks and to generalize more efficiently under data-limited conditions. 

\begin{table}[t]
\centering
\small
\setlength{\tabcolsep}{6pt}
\begin{tabular}{ccc}
\toprule
\textbf{Model} & \textbf{Training Time (s)} $[10^{2}]$ & \textbf{Prediction Time (s) $[10^{-2}]$} \\
\midrule
GP $s=1$ & $6.38$  & 4.0 \\
GP $s=2$ & $11.40$ & 3.9 \\
GP $s=3$ & $7.20$  & 4.0 \\
GP $s=4$ & $21.50$ & 4.3 \\
\midrule
\textbf{MTGP} & $1.66$ & 50.3 \\
\bottomrule
\end{tabular}
\caption{Average training and prediction times per test scenario for the single-task GP and the MTGP models.}
\label{tab:gp_mtgp_times}
\end{table}

\begin{figure}[t!]
    \centering
    \includegraphics[width=0.55\textwidth]{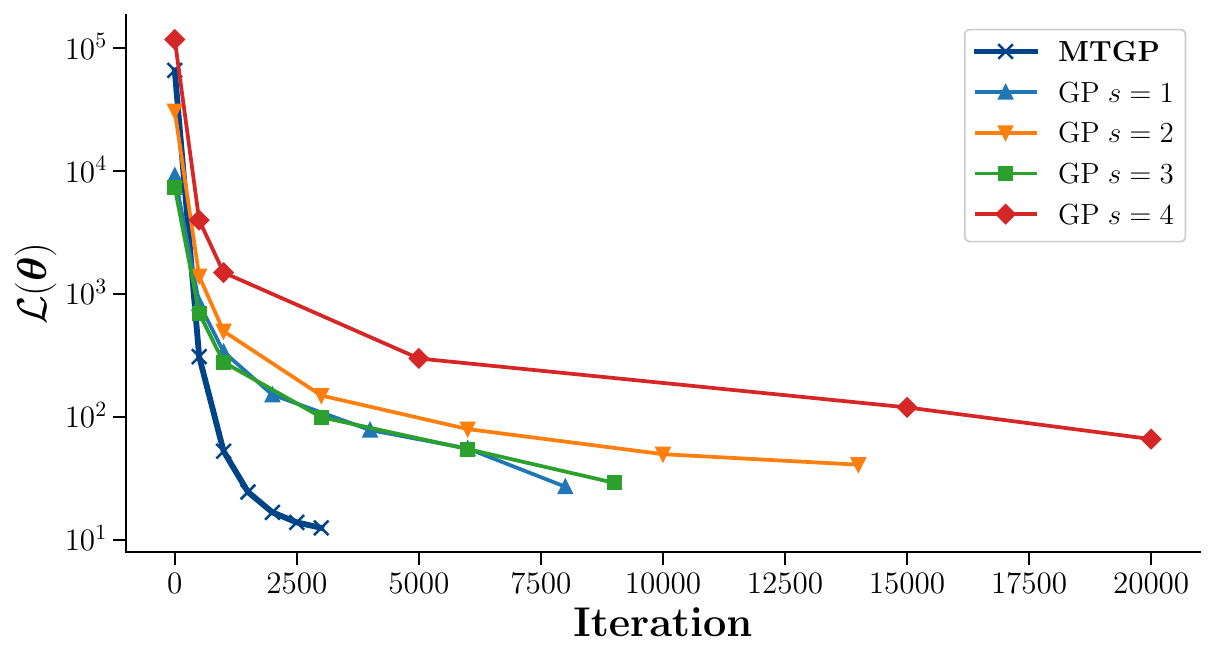}
    \caption{Optimization convergence comparison between independently trained single-task GP models and the MTGP. The reported curves correspond to the average negative log-marginal likelihood $\mathcal{L}(\boldsymbol{\theta})$ computed over five  runs, plotted as a function of the number of training iterations.}
    \label{fig:convergence_mono_vs_mtgp}
\end{figure}

From a computational standpoint, we observe a reduction in training time until optimization convergence as reported in Table~\ref{tab:gp_mtgp_times}. While the CPU times for the four independently trained GPs range from approximately 10 minutes to more than 30 minutes, the one for MTGP is about 3 minutes on average. All reported times correspond to the mean over five independent runs with different random initialization, with variations of approximately $1\%$ in both the final training time and the attained log-marginal likelihood values, indicating reliable repeatability of the optimization procedure. This gain in computing efficiency was not anticipated as the MTGP model has more parameters to learn than the GP. 

Such empirical gain is explained by the convergence behavior shown in Figure~\ref{fig:convergence_mono_vs_mtgp}, where we observe that the MTGP rapidly decreases the negative log-marginal likelihood and reaches a near-stationary regime after approximately $3\times 10^{3}$ iterations. In contrast, the single-task GP models exhibit a slower decay, with some tasks requiring more than $10^{4}$ iterations to reach comparable log-likelihood values. This behavior can be interpreted in light of the theoretical properties of multitask learning curves. In the presence of non-zero inter-task correlations, the prediction error on a given task depends not only on its own sample size but also on the number of observations available for correlated tasks. When tasks are moderately to strongly correlated, an initial collective learning phase occurs, during which the error decreases rapidly due to shared statistical structure. This translates into a more favorable marginal likelihood landscape for the joint model, yielding faster practical convergence. This is indeed the case in our study where tasks $y_2$, $y_3$, and $y_4$ exhibit non-negligible positive and negative correlations. We can then conclude that the acceleration observed here is due to the explicit coupling of tasks through the learned inter-task covariance matrix $\mathbf{K}_{\mathcal S}$. Its off-diagonal entries quantify statistical transfer, and its spectral structure determines the dominant shared modes across tasks. Joint hyperparameter estimation concentrates learning in a lower-dimensional shared subspace, whereas single-task GPs repeatedly estimate similar structural patterns separately, leading to redundant optimization effort.

The situation differs at the prediction stage. As reported in
Table~\ref{tab:gp_mtgp_times}, the average prediction time of the
single-task GPs is $0.04\,\mathrm{s}$, whereas the MTGP requires
$0.50\,\mathrm{s}$ to jointly predict the four outputs. For a single-task GP, prediction requires solving a linear system
associated with $n_f n_u = 8112$ degrees of freedom for each task independently. In contrast, the MTGP jointly models the
$S=4$ correlated tasks, leading to a global effective system size of $Sn_fn_u=32448$. This fourfold increase in effective dimension directly explains the higher computational cost at prediction time. Although the Kronecker factorization ensures that computations remain tractable and memory-efficient, the coherent coupling of all four outputs within a unified multitask framework naturally entails a higher computational burden than fully decoupled single-task models. More precisely, prediction still requires coupled
linear solves and structured matrix--vector operations involving all
tasks simultaneously. This interpretation was confirmed by the CPU profiling results.
For a single MTGP prediction, more than $50\%$ of the computation time
is spent in dense matrix multiplications, and approximately $44\%$ in
linear system solves, while kernel evaluation itself represents a
negligible fraction of the total cost.

\section{Conclusions and future work}
\label{sec:conclusion}

This work introduces a scalable multitask Gaussian process (MTGP) designed for systems whose behavior depends on functional covariates.
In contrast with single-task GPs, which ignore structural correlations between outputs, the multitask formulation explicitly exploits cross-task dependencies through a separable kernel.
These assumptions result in a Kronecker-structured covariance that has been exploited for MTGP computational efficiency. The combined use of MTGPs and functional covariates offers a general approach for representing variability effects in mechanical systems, particularly when experimental data are costly and the responses exhibit complex, multi-modal dynamics. Moreover, the specific treatment of the Kronecker covariance structure that we have proposed allows a significant computational gain, making exact GP inference feasible at scales that would be prohibitively expensive for dense models.

The effectiveness of the MTGP has been demonstrated on a synthetic benchmark and on a real mechanical riveted assembly. 
In this application, the functional covariates encode rivets material variability through displacement-dependent characteristics, while the outputs correspond to force-displacement relations at several locations in the assembly. The MTGP can accurately learn from 50 to 80 samples the complex dynamic behavior of the assembly, including oscillatory patterns and task-dependent amplitude variations. It also provides calibrated confidence intervals in less than a CPU second on a current standard computer, enabling design exploration and uncertainty propagation in such mechanical assemblies structures.

With regard to model learning, the experiments carried out have shown two phenomena. It was observed that, for the MTGP model, predicting an accurate mean response takes fewer data than predicting a representative uncertainty interval. Also, even though the MTGP model has more parameters than single GPs, learning them may actually be faster thanks a more favorable likelihood landscape.

Future work may consider extensions to multivariate functional covariates, such as those representing spatial fields (e.g., surface topographies). While such inputs would allow a richer description of mechanical variabilities, they also introduce substantially higher computational and modeling complexity. Therefore, to maintain a tractable inference, it is required to develop even more scalable MTGP models that benefit, in addition to a structured covariance, from low-dimensional representations of the multivariate functions.

\section*{Acknowledgements}
This research has been supported by the project GAME (ANR-23-CE46-0007), funded by the French National Research Agency (ANR).
Part of the work has been conducted when AFLL was affiliated at CERAMATHS, UPHF. The authors gratefully acknowledge Dr. Nicolas Leconte (Researcher at the French Aerospace Lab, ONERA, DMAS, Lille) for providing the numerical model used in this study.

\newpage
\appendix

\section{Description of the encoding strategies used for functional inputs}
\label{app:encodings}
In accordance with the projection framework of Section~\ref{subsec:approximation-functional}, we present the four encoding strategies used for functional inputs in this work, organized from direct data-driven approaches to basis-based and hybrid representations.

\subsection{PCA}
For each functional component $f_d \in L^2(\mathcal{T})$, with $d \in \{1,\dots,d_f\}$, 
the function is sampled on a common grid 
$\{u_1,\dots,u_{n_u}\} \subset \mathcal{T}$, yielding for each replicate 
$i \in \{1,\dots,n_f\}$ the discretized vector
$
\bm{f}_{d,i} = \begin{bmatrix}f_{d,i}(u_1),\dots,f_{d,i}(u_{n_u})\end{bmatrix}^\top \in \mathbb{R}^{n_u}$.

Stacking all replicates gives the centered data matrix
$
\bm{\mathcal{F}}_d =
\begin{bmatrix}
\bm{f}_{d,1} & \cdots & \bm{f}_{d,n_f}
\end{bmatrix}^\top
\in \mathbb{R}^{n_f \times n_u}$. The empirical covariance operator, expressed in discretized form, is given by
\begin{equation}
	\bm{\kappa}_d = \frac{1}{n_f}\,\bm{\mathcal{F}}_d^\top \bm{\mathcal{F}}_d
	\in \mathbb{R}^{n_u \times n_u}.
\end{equation}
Let $\{\Upsilon_{d,r}\}_{r=1}^{n_u}$ denote the orthonormal eigenfunctions of this operator, 
identified through their evaluations on the grid $\{u_k\}_{k=1}^{n_u}$, with associated 
eigenvalues $\{\lambda_{d,r}\}_{r=1}^{n_u}$ sorted in decreasing order.

Let $p_d^\ast$ be the smallest integer such that
\begin{equation}
	\sum_{r=1}^{p_d^\ast}
	\frac{\lambda_{d,r}}{\sum_{r'=1}^{n_u} \lambda_{d,r'}}
	\;\geq\; \mathcal{I}_\ast,
\end{equation}
where $\mathcal{I}_\ast \in (0,1)$ is a prescribed inertia threshold.
The PCA approximation of $f_{d,i}$ is then written in the form
\begin{equation}
	f_{d,i}(u) \;\approx\;
	\sum_{r=1}^{p_d^\ast} \beta_{d,r}^{(i)}\,\Upsilon_{d,r}(u),
\end{equation}
where the coefficients
$
\beta_{d,r}^{(i)} =
\int_{\mathcal{T}} f_{d,i}(u)\,\Upsilon_{d,r}(u) \ du$ are approximated numerically from the discretized observations. The vector $ \bm{\beta}_d = [\beta_{d,1}, \ldots, \beta_{d,p_d}]^\top $
thus provides a data-driven low-dimensional surrogate representation of $f_d$, consistent 
with the general projection framework introduced in 
Section~\ref{subsec:approximation-functional}.

\subsection{Projection on functional bases}
As detailed in Section~\ref{subsec:approximation-functional},
each $f_d \in L^2(\mathcal{T})$ is approximated in a finite-dimensional subspace
spanned by a chosen basis $\{\Upsilon_{d,r}\}_{r=1}^{p_d}$, yielding a coefficient
vector $\bm{\beta}_d \in \mathbb{R}^{p_d}$.

\subsubsection{B-spline basis}

In this spline-based setting, the functional input $f_d$ is approximated through
the general projection
\[
f_d(u) \approx \sum_{r=1}^{p_d} \beta_{d,r}\,\Upsilon_{d,r}(u),
\]
where $\{\Upsilon_{d,r}\}_{r=1}^{p_d}$ denotes a B-spline basis of order $m$
(i.e.\ piecewise polynomials of degree $m-1$).

The basis functions are defined over a non-decreasing knot sequence
$\{\tau_{d,r}\}_{r=1}^{p_d+m}$ satisfying
\[
\tau_{d,1} = \cdots = \tau_{d,m}
< \tau_{d,m+1} < \cdots < \tau_{d,p_d}
< \tau_{d,p_d+1} = \cdots = \tau_{d,p_d+m}.
\]

This construction corresponds to a clamped knot vector, involving
$p_d - m$ distinct interior knots and $m$ repeated knots at each boundary.
Under the assumption of simple interior knots, the resulting B-spline basis
is globally $C^{m-2}$-continuous over $\mathcal{T}$.

The basis functions are defined recursively via the Cox--de Boor formula.
For order~$1$,
\[
\Upsilon_{d,r}^{(1)}(u)
= \mathbf{1}_{[\tau_{d,r},\,\tau_{d,r+1})}(u),
\qquad r \in \{1,\dots,p_d+m-1\},
\]
and for orders $m \geq 2$,
\[
\Upsilon_{d,r}^{(m)}(u)
=
\frac{u-\tau_{d,r}}{\tau_{d,r+m-1}-\tau_{d,r}}\,\Upsilon_{d,r}^{(m-1)}(u)
+
\frac{\tau_{d,r+m}-u}{\tau_{d,r+m}-\tau_{d,r+1}}\,\Upsilon_{d,r+1}^{(m-1)}(u),
\qquad r \in \{1,\dots,p_d\},
\]
with the convention that $\Upsilon_{d,r}^{(m)}(u)=0$ whenever a denominator
vanishes.

\medskip
The resulting basis functions $\{\Upsilon_{d,r}\}_{r=1}^{p_d}$ are nonnegative,
compactly supported, and form a partition of unity:
\[
\sum_{r=1}^{p_d} \Upsilon_{d,r}(u) = 1,
\quad \forall\, u \in \mathcal{T}.
\]
These properties make B-spline bases particularly suitable for the approximation
of smooth functional inputs, as they provide local control, numerical stability,
and global smoothness of order $C^{m-2}$.

\subsubsection{Wavelet basis}

Wavelet-based representations rely on a multiresolution decomposition of
$f_d \in L^2(\mathcal{T})$.
A wavelet system is defined by a scaling function $\varphi_d$ and a mother wavelet
$\psi_d$, which generate families of basis functions through dyadic dilations and
translations.
To unify notation, we introduce the generic wavelet atom
\[
\Upsilon_{d,j,k}(u)
=
2^{j/2}\,\Upsilon_d(2^j u - k),
\qquad j,k \in \mathbb{Z},
\]
where $\Upsilon_d$ denotes either the scaling function $\varphi_d$ (approximation atoms)
or the mother wavelet $\psi_d$ (detail atoms).
Here, $j$ denotes the scale index and $k$ the translation index.
Since the functional domain $\mathcal{T}$ is bounded, only a finite number of
translations are retained at each scale; more precisely, for each $j$, the index
$k$ ranges over a finite subset $\mathcal{K}_j \subset \mathbb{Z}$ such that the support
of the corresponding atom $\Upsilon_{d,j,k}$ intersects $\mathcal{T}$.

A truncated multiresolution expansion is obtained by retaining the approximation space
at level $j_0$ and the detail spaces up to level $j_{\max}$, leading to the representation

\begin{equation}
	f_d(u) \;\simeq\;
	\sum_{k \in \mathcal{K}_{j_0}} c_{d,j_0,k}\,\varphi_{d,j_0,k}(u)
	+
	\sum_{j=j_0}^{j_{\max}} \sum_{k \in \mathcal{K}_j}
	d_{d,j,k}\,\psi_{d,j,k}(u).
	\label{eq:wavelet-truncated}
\end{equation}

The associated wavelet atoms
$\{\varphi_{d,j_0,k}\}_{k \in \mathcal{K}_{j_0}} \cup
\{\psi_{d,j,k}\}_{j=j_0,\dots,j_{\max},\,k \in \mathcal{K}_j}= \{\Upsilon_{d,r}\}_{r=1}^{p_d}$
constitute a finite basis, whose expansion coefficients
are collected in the vector $\boldsymbol{\beta}_d \in \mathbb{R}^{p_d}$.

For orthonormal wavelet systems (e.g., Haar or Daubechies), the basis
$\{\Upsilon_{d,r}\}_{r=1}^{p_d}$ is orthonormal in $L^2(\mathcal{T})$.
As a consequence, the $L^2$ distance between two functions reduces to the Euclidean
distance between their corresponding coefficient vectors.

\subsection{PCA on projection coefficients}

After projection onto a fixed functional basis
$\{\Upsilon_{d,r}\}_{r=1}^{p_d}$ (e.g., spline or wavelet bases),
each replicate of the functional component $f_d$ is represented by its coefficient
vector $\boldsymbol{\beta}_{d,i} \in \mathbb{R}^{p_d}$, for
$i \in \{1,\dots,n_f\}$.
Let $\{\boldsymbol{\beta}_{d,i}\}_{i=1}^{n_f}$ denote the resulting sample of
projection coefficients, centered componentwise.

PCA is then applied to this collection of vectors in order
to identify a low-dimensional subspace capturing the dominant modes of variability.
Let $p_d^\ast < p_d$ be the smallest integer satisfying the inertia criterion
associated with a prescribed threshold $\mathcal{I}_\ast$.
Each replicate is subsequently encoded by a reduced vector
$\boldsymbol{\gamma}_{d,i} \in \mathbb{R}^{p_d^\ast}$, obtained by projection of
$\boldsymbol{\beta}_{d,i}$ onto the leading principal directions.

This two-step procedure functional projection followed by PCA on the resulting
coefficients yields a compact, decorrelated, and basis-aware representation of the
functional input, while preserving the most informative directions induced by the
chosen functional basis.

\section{Complementary numerical results: prediction results for the different functional encodings}
\label{app:comp_results}

\begin{figure}[t!]
    \centering
 \begin{subfigure}{\textwidth}
        \centering
        \rotatebox{90}{\footnotesize \hspace{-3ex}B-spline}
        \vspace{0.3em}
        \begin{minipage}{0.32\textwidth}
            \centering
            \includegraphics[width=\linewidth]{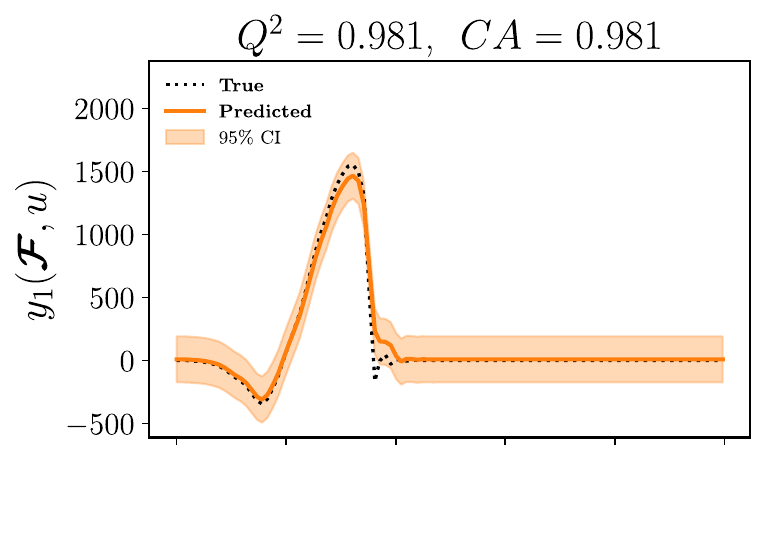}
        \end{minipage}%
        \begin{minipage}{0.32\textwidth}
            \centering
            \includegraphics[width=\linewidth]{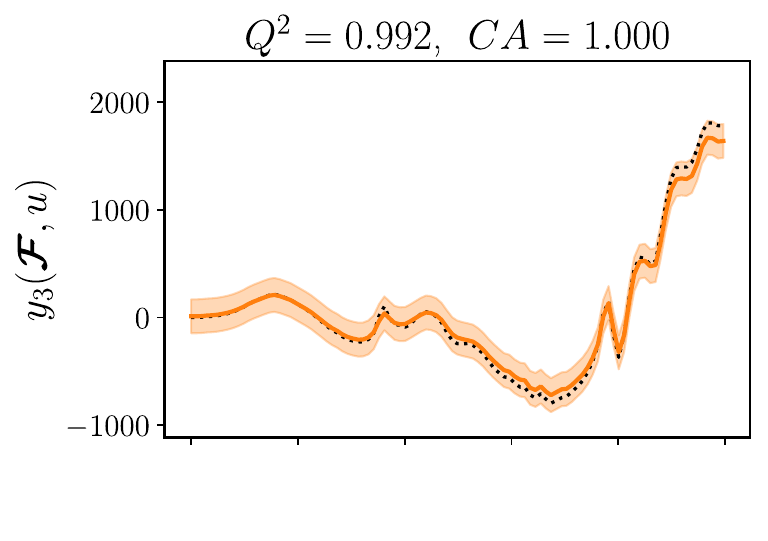}
        \end{minipage}%
        \begin{minipage}{0.32\textwidth}
            \centering
            \includegraphics[width=\linewidth]{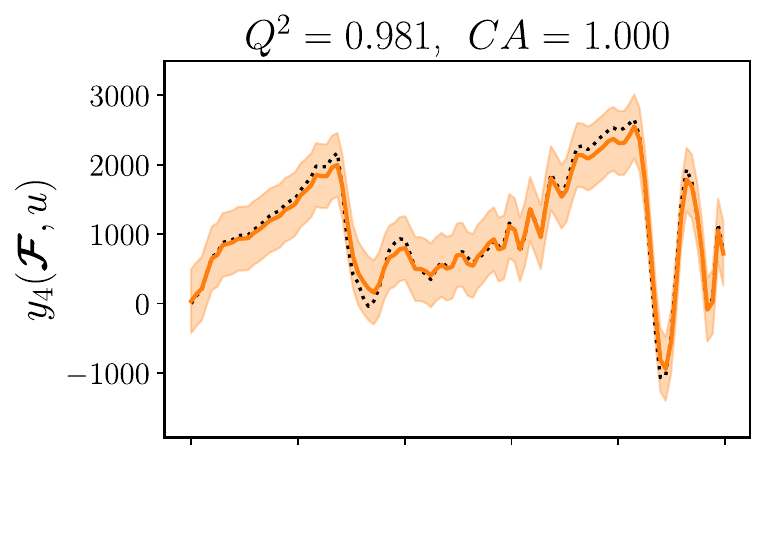}
        \end{minipage}
        \vspace{-0.5em}
    \end{subfigure}

    \begin{subfigure}{\textwidth}
        \centering
        \rotatebox{90}{\footnotesize \hspace{-7ex}B-spline + PCA}
        \vspace{0.3em}
        \begin{minipage}{0.32\textwidth}
            \centering
            \includegraphics[width=\linewidth]{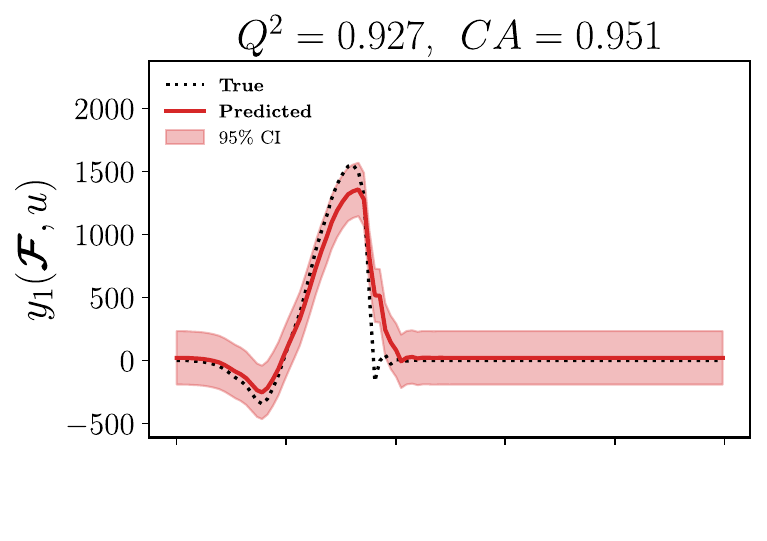}
        \end{minipage}%
        \begin{minipage}{0.32\textwidth}
            \centering
            \includegraphics[width=\linewidth]{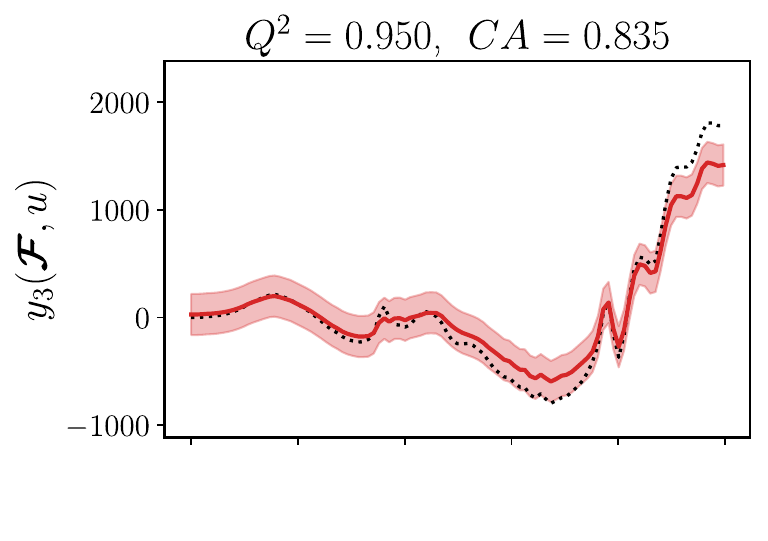}
        \end{minipage}%
        \begin{minipage}{0.32\textwidth}
            \centering
            \includegraphics[width=\linewidth]{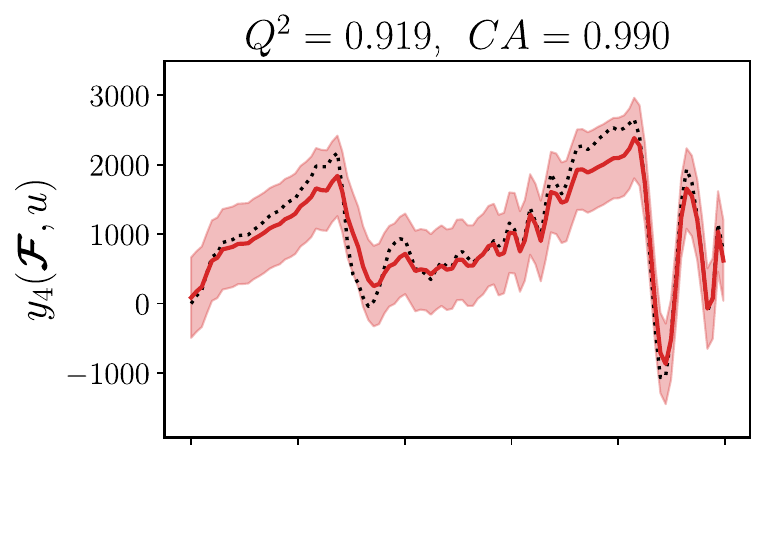}
        \end{minipage}
        \vspace{-0.5em}
    \end{subfigure}

    \begin{subfigure}{\textwidth}
        \centering
        \rotatebox{90}{\footnotesize \hspace{-3ex}Wavelets}
       \vspace{0.3em}
        \begin{minipage}{0.32\textwidth}
            \centering
            \includegraphics[width=\linewidth]{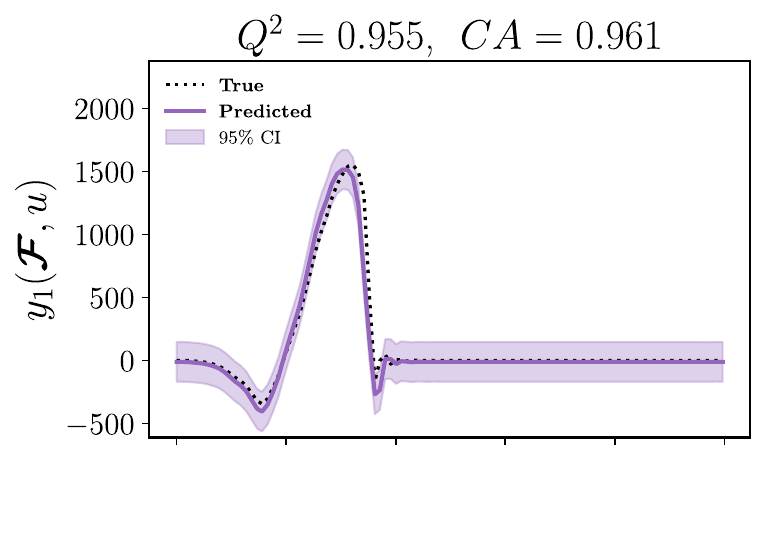}
        \end{minipage}%
        \begin{minipage}{0.32\textwidth}
            \centering
            \includegraphics[width=\linewidth]{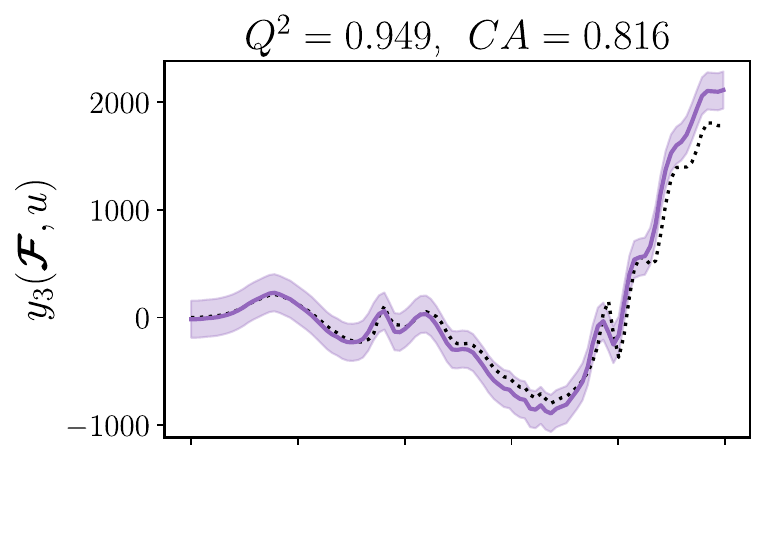}
        \end{minipage}%
        \begin{minipage}{0.32\textwidth}
            \centering
            \includegraphics[width=\linewidth]{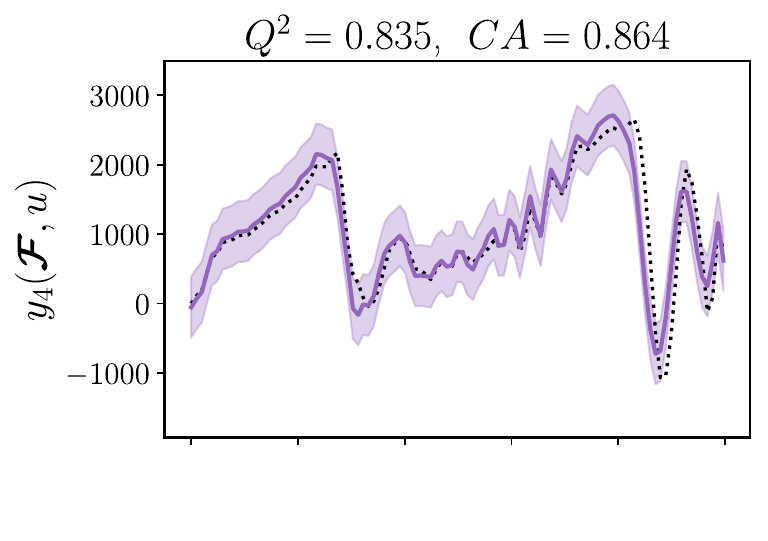}
        \end{minipage}
        \vspace{-0.5em}
    \end{subfigure}

    \begin{subfigure}{\textwidth}
        \centering
        \rotatebox{90}{\footnotesize \hspace{-7.5ex}Wavelets + PCA}
        \vspace{0.3em}
        \begin{minipage}{0.32\textwidth}
            \centering
            \includegraphics[width=\linewidth]{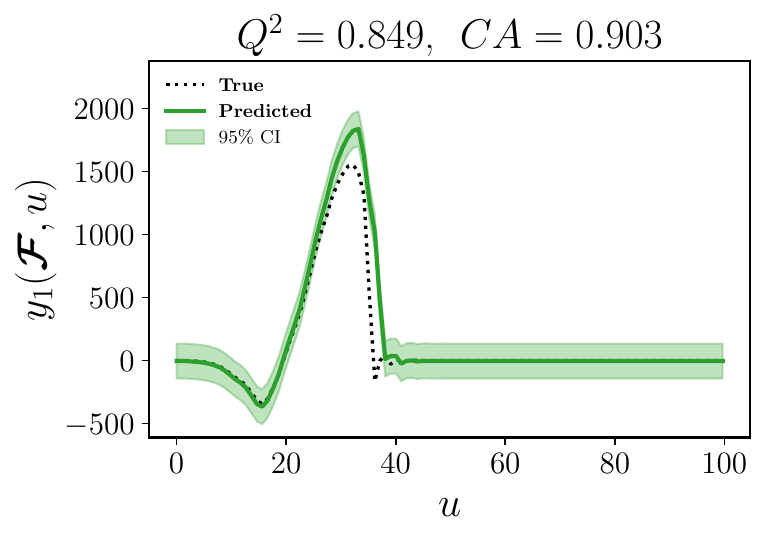}
        \end{minipage}%
        \begin{minipage}{0.32\textwidth}
            \centering
            \includegraphics[width=\linewidth]{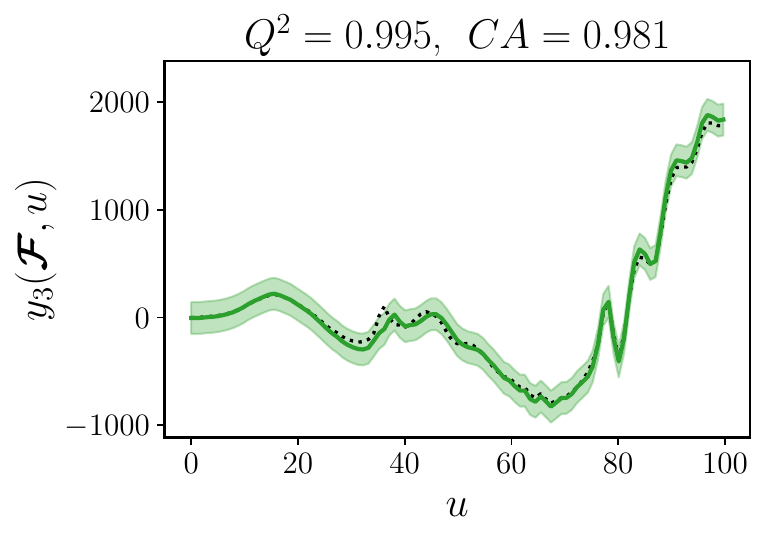}
        \end{minipage}%
        \begin{minipage}{0.32\textwidth}
            \centering
            \includegraphics[width=\linewidth]{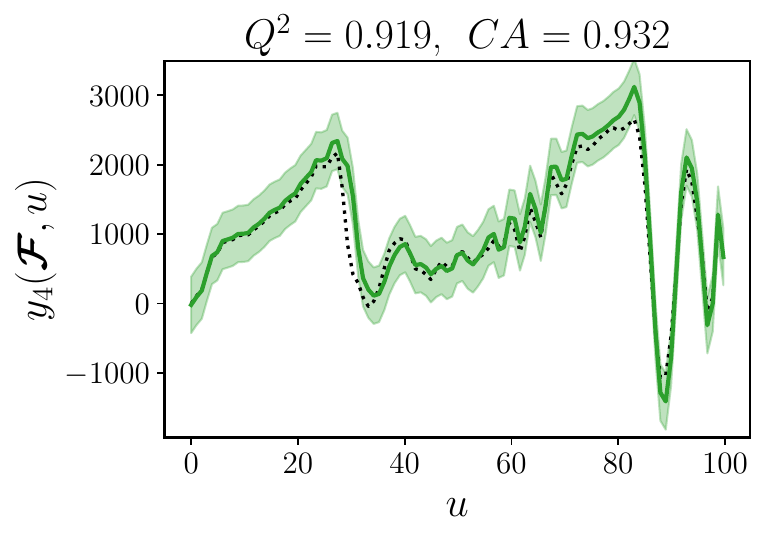}
        \end{minipage}
    \end{subfigure}

    \caption{
        MTGP predictions for Scenario~23 across the tasks $s \in \{1, 3, 4\}$ (indexed by columns) and the four functional encodings considered in Section~\ref{subsec:gp_perf_training_size} (indexed by rows).} 
    \label{fig:fmtgp_appendix_pred}
\end{figure}
Figure~\ref{fig:fmtgp_appendix_pred} compares the predictions obtained with four functional
encoding strategies B-splines, B-splines + PCA, wavelets, and wavelets + PCA for
Scenario~23 with \(n_f = 78\) training functions. The results reveal contrasted behaviors across tasks and encodings.
B-spline representations without PCA provide overall strong predictive performance,
particularly for tasks \(s=3\) and \(s=4\), with high \(Q^2\) values and satisfactory
coverage.
By contrast, applying PCA to B-spline coefficients leads to a visible degradation of
performance for several tasks, including a marked decrease in both \(Q^2\) and coverage
for \(s=3\). Wavelet-based representations without PCA exhibit larger task-dependent variability,
with reasonable accuracy for \(s=1\) and \(s=2\), but a noticeable loss of precision for
\(s=4\).
Applying PCA to wavelet coefficients improves predictions for certain tasks,
most notably \(y_3\), where both the coefficient of determination and coverage increase,
while improvements remain limited or negligible for other tasks.

These observations indicate that the impact of PCA strongly depends on both the functional
representation and the output considered. In particular, PCA can be beneficial for
wavelet-based encodings in some cases, but may be counterproductive when applied to
already well-structured representations such as B-splines.

\section{Efficient kronecker-based computations and mode-wise vectorization}

\subsection{Computation of \texorpdfstring{$\log |\mathbf{L}|$}{det(L)}}
\label{appendix:Kronecker:logL}

Using the standard determinant identity for Kronecker products $|\mathbf{A} \otimes \mathbf{B}| = |\mathbf{A}|^{\dim(B)} \, |\mathbf{B}|^{\dim(A)}$, and applying it recursively, we obtain a closed-form expression for the determinant of $\mathbf{L} = \mathbf{L}_{\mathcal{S}} \otimes \mathbf{L}_{f} \otimes \mathbf{L}_{u}$:
\begin{align*}
    |\mathbf{L}|
    &= 
    |\mathbf{L}_{\mathcal{S}}|^{\,n_f n_u}
    \; |\mathbf{L}_{f}|^{\,S n_u}
    \; |\mathbf{L}_{u}|^{\,S n_f} \\
    &= 
    \left( \prod_{i=1}^{S} (\mathbf{L}_{\mathcal{S}})_{ii} \right)^{n_f n_u}
    \left( \prod_{i=1}^{n_f} (\mathbf{L}_{f})_{ii} \right)^{S n_u}
    \left( \prod_{i=1}^{n_u} (\mathbf{L}_{u})_{ii} \right)^{S n_f}.
\end{align*}
Then, by applying the property $\log(ab)^c = c (\log a + \log b)$, we finally obtain
\begin{equation*}
    \log |\mathbf{L}|
    = (n_f n_u) \sum_{i=1}^{S} \log (\mathbf{L}_{\mathcal{S}})_{ii}
    + (S n_u) \sum_{i=1}^{n_f} \log (\mathbf{L}_{f})_{ii}
    + (S n_f) \sum_{i=1}^{n_u} \log (\mathbf{L}_{u})_{ii}.
\end{equation*}

\subsection{Kronecker--vec identity.}
\label{Kronecker--vec}
We adopt the standard column-major vectorization. 
For $\mathcal{Y} \in \mathbb{R}^{n_1 \times \cdots \times n_D}$ 
and matrices $A_d \in \mathbb{R}^{m_d \times n_d}$, the following identity holds:
\begin{equation}
	\left( \bigotimes_{d=1}^{D} A_d \right)\, \operatorname{vec}(\mathcal{Y})
	= \operatorname{vec}\!\Big( \mathcal{Y} \times_D A_D \times_{D-1} A_{D-1} \cdots \times_1 A_1 \Big),
	\label{eq:KronProd}
\end{equation}
where $\times_d$ denotes the mode-$d$ product. 

\begin{itemize}
    \item \textbf{Case $D=2$.}
Let $Y \in \mathbb{R}^{n_1 \times n_2}$, 
$A_1 \in \mathbb{R}^{m_1 \times n_1}$, 
$A_2 \in \mathbb{R}^{m_2 \times n_2}$. Then
\[
(A_1 \otimes A_2)\, \operatorname{vec}(Y)
= \operatorname{vec}\big(A_2\, Y\, A_1^\top\big).
\]
This shows that with the ordering $\bigotimes_{d=1}^2 A_d$, 
the action of $A_1$ appears on the \emph{last} mode.

\item \textbf{Case $D=3$.}
Let $\mathcal{Y} \in \mathbb{R}^{n_1 \times n_2 \times n_3}$ 
and $A_d \in \mathbb{R}^{m_d \times n_d}$. Then
\[
\operatorname{vec}(\mathcal{Y} \times_3 A_3 \times_2 A_2 \times_1 A_1)
= (A_1 \otimes A_2 \otimes A_3)\, \operatorname{vec}(\mathcal{Y}).
\]
Step-by-step application gives
\begin{align*}
	\operatorname{vec}(\mathcal{Y} \times_3 A_3) 
	&= (A_3 \otimes I_{n_2} \otimes I_{n_1})\, \operatorname{vec}(\mathcal{Y}),\\
	\operatorname{vec}(\mathcal{Y} \times_3 A_3 \times_2 A_2) 
	&= (A_3 \otimes A_2 \otimes I_{n_1})\, \operatorname{vec}(\mathcal{Y}),\\
	\operatorname{vec}(\mathcal{Y} \times_3 A_3 \times_2 A_2 \times_1 A_1) 
	&= (A_1 \otimes A_2 \otimes A_3)\, \operatorname{vec}(\mathcal{Y}).
\end{align*}

\item \textbf{General case $D > 3$ by induction.}
Assume~\eqref{eq:KronProd} holds for $D-1$. Then
\begin{align*}
	\operatorname{vec}\!\Big( (\mathcal{Y} \times_{D-1} A_{D-1} \cdots \times_1 A_1) \times_D A_D \Big)
	&= (A_D \otimes I)\, \operatorname{vec}\!\Big( \mathcal{Y} \times_{D-1} A_{D-1} \cdots \times_1 A_1 \Big) \\
	&= (A_D \otimes I)\, \left( \bigotimes_{d=1}^{D-1} A_d \right)\, \operatorname{vec}(\mathcal{Y}) \\
	&= \left( \bigotimes_{d=1}^{D} A_d \right)\, \operatorname{vec}(\mathcal{Y}).
\end{align*}
\end{itemize}

\medskip

\noindent\textbf{Remark.}
With the convention $\bigotimes_{d=1}^D A_d$, the Kronecker factors appear
in increasing order on the left of~\eqref{eq:KronProd}, while the mode-wise products appear in decreasing order on the right. This explains the apparent reversal. For clarity, we adopt this convention consistently throughout the manuscript.

\bibliographystyle{elsarticle-harv}
\bibliography{references}

\end{document}